\input ./style/arxiv-general.cfg
\documentclass[aos,MSNbibl,nameyear,seceqn,dvips]{arximspdf}
\makeatletter
   \@ifpackageloaded{graphicx}{}{\usepackage{graphicx}}
\makeatother
\usepackage{amssymb,eufrak,euscript}

\doi{10.1214/15-AOS1315}
\volume{43}
\issue{4}
\pubyear{2015}
\firstpage{1535}
\lastpage{1567}
\docsubty{FLA}

\makeatletter

\def\vfrac#1#2{(#1)/#2}

\def\sklfrac#1#2{(#1/#2)}

\newcommand{\lleft}{\left}
\newcommand{\rrvert}{\vert}
\newcommand{\rright}{\right}
\newcommand{\rrVert}{\Vert}
\newcommand{\llvert}{\vert}
\newcommand{\llVert}{\Vert}
\newtheorem{teo}{Theorem}[section]
\newproclaim{assumption}[teo]{Assumption}
\newtheorem{prop}[teo]{Proposition}
\newproclaim{exa}{Example}
\newproclaim{rems}{Remarks}
\newproclaim{rem}{Remark}
\makeatother

\begin{document}
\begin{frontmatter}

\title{Regularized estimation in sparse high-dimensional time series models}
\runtitle{Regularized estimation in time series}

\begin{aug}
\author[A]{\fnms{Sumanta}~\snm{Basu}\ead[label=e1]{sumbose@umich.edu}} 
\and
\author[A]{\fnms{George}~\snm{Michailidis}\corref{}\ead[label=e2]{gmichail@umich.edu}\thanksref{T1}}
\runauthor{S. Basu and G. Michailidis}
\affiliation{University of Michigan}
\address[A]{Department of Statistics\\
University of Michigan\\
Ann Arbor Michigan 48109\\
USA\\
\printead{e1}\\
\phantom{E-mail: }\printead*{e2}}
\end{aug}
\thankstext{T1}{Supported by NSA Grant H98230-10-1-0203 and NSF Grants DMS-11-61838 and DMS-12-28164.}
\runauthor{Sumanta Basu and George Michailidis}

%
\received{\smonth{2} \syear{2014}}
%
\revised{\smonth{1} \syear{2015}}

%
\begin{abstract}
Many scientific and economic problems involve the analysis of
high-dimensional time series datasets.
However, theoretical studies in high-dimensional statistics to date
rely primarily on the assumption of
\mbox{independent} and identically distributed (i.i.d.) samples. In this work,
we focus on stable Gaussian
processes and investigate the theoretical properties of \mbox{$\ell
_1$-}regularized estimates in two important
statistical problems in the context of high-dimensional time series:
(a) stochastic regression with
serially correlated errors and (b) transition matrix estimation in
vector autoregressive (VAR) models.
We derive nonasymptotic upper bounds on the estimation errors of the
regularized estimates and
establish that consistent estimation under high-dimensional scaling is
possible via {$\ell_1$-}regularization
for a large class of stable processes under sparsity constraints. A key
technical contribution of the work is
to introduce a measure of stability for stationary processes using
their spectral properties that provides insight into the effect of
dependence on the accuracy of the regularized estimates. With this
proposed stability measure, we establish some useful deviation bounds
for dependent data, which can be used to study several important
regularized estimates in a time series setting.
\end{abstract}

%
\end{frontmatter}

\section{Introduction}\label{secintro}

Recent advances in information technology have made high-dimensional
time series data sets increasingly common in numerous applications.
Examples include structural analysis and forecasting with a large
number of macroeconomic variables [\citet{demol08}],
reconstruction of gene regulatory networks from time course microarray
data [\citet{Michailidis-review}], portfolio selection and
volatility matrix estimation in finance [\citet{fan-lv-08}] and
studying co-activation networks in human brains using task-based or
resting state fMRI data [\citet{Smith-12-neuroimage}]. These
applications require analyzing a large number of temporally observed
variables using small to moderate sample sizes (number of time points),
and the techniques used for the respective learning tasks
include classical regression, vector autorgressive modeling and
covariance estimation. Meaningful inference in such settings is often
impossible without imposing some lower-dimensional structural
assumption on the data generating mechanism, the most common being that
of sparsity on the model parameter space. In high-dimensional
regression and VAR problems, the notion of sparsity is often
incorporated into the estimation procedure by \mbox{$\ell_1$-}penalization
procedures like lasso and its variants [\citet
{bickel2009simultaneous,vandegeerthreshadaptejs2011}], while for
covariance matrix estimation problems, sparsity is enforced via hard
thresholding [\citet{BL08covthresh}].

Theoretical properties of such regularized estimates under
high-\break dimensional scaling have been investigated in numerous studies
over the last few years, under the key assumption that the samples are
independent and identically distributed (i.i.d.). On the other hand,
theoretical analysis of these estimates in a time series context, where
the data exhibit \textit{temporal} and \textit{cross-sectional dependence},
is rather incomplete. A central challenge is to assess how the
underlying dependence structure affects the performance of these
regularized estimates.

In this paper, we focus on stationary Gaussian time series and use
their \textit{spectral properties} to propose a measure of stability.
Using this measure of stability, we establish necessary concentration
bounds for dependent data and study, in a nonasymptotic framework, the
theoretical properties of regularized estimates in the following key
statistical models: (a) $\ell_1$-penalized sparse stochastic
regression with exogenous predictors and serially correlated errors and
(b) $\ell_1$-penalized least squares and log likelihood based
estimation of sparse VAR models. We establish nonasymptotic upper
bounds on the estimation error and show that lasso can perform
consistent estimation in high-dimensional settings under a mild
stability assumption on the underlying processes that is common in the
classical literature of low-dimensional time series. Our results also
provide new insights into how the convergence rates are affected by the
presence of temporal dependence in the data.

Next, we introduce the two models analyzed in this paper and highlight
the main contributions of our work to the existing literature. Although
the main interest of this work is to study VAR models in
high dimensions, a key stepping stone to our analysis comes from
stochastic regression models, which are of independent interest.

\subsection*{Stochastic regression}
We start with this canonical problem
in time series analysis [\citet{hamilton1994time}],
a linear regression model of the form
%
\begin{equation}
\label{eqnstoch-reg-defn-intro} y^t = \bigl\langle\beta^*, X^t \bigr
\rangle+ \varepsilon^t,\qquad t = 1, \ldots, n,
\end{equation}
where the $p$-dimensional predictors $\{X^t \}$ and the errors $\{
\varepsilon^t \}$ are generated according to independent,
centered, Gaussian stationary processes. Under a sparsity assumption on
$\beta^*$, we study the properties of the lasso estimate
%
\begin{equation}
\label{eqnlasso-intro} \hat{\beta} = \mathop{\operatorname{argmin}}_{\beta\in\mathbb
{R}^p}
\frac{1}{n} \llVert Y - \mathcal{X} \beta\rrVert ^2 +
\lambda_n \llVert \beta \rrVert _1,
\end{equation}
where $Y = [y^n: \ldots: y^1]'$, $\mathcal{X} = [X^n: \ldots: X^1]'$
and $\llVert   \beta\rrVert  _1 = \sum_{j = 1}^p \llvert  \beta_j\rrvert  $.
Theoretical properties of lasso have been studied for fixed design
regression $Y = \mathcal{X} \beta^* + E$, with $E = [e^n: \ldots:
e^1]'$, by several authors [\citet{bickel2009simultaneous,powai2012,negahban2012unified}]. They establish consistency of lasso
estimates in a high-dimensional regime under some form of restricted
eigenvalue (RE) or restricted strong convexity (RSC) assumption on $S =
\mathcal{X}'\mathcal{X}/n$ and suitable deviation conditions on
$\mathcal{X}'E/n$.

In general, for a given design matrix $\mathcal{X}$, verifying that
$\mathcal{X}$ satisfies an RE condition [\citet
{dobriban2013regularity}] is an NP-hard problem. In the case where the
rows of $\mathcal{X}$ are independently generated from a common
Gaussian/sub-Gaussian ensemble, these assumptions are known to hold
with high probability under mild conditions [\citet
{raskutti2010REcorrgauss,rudelsonzhou2012paper}]. It is not clear,
however, whether similar regularity conditions are satisfied with high
probability when the observations are dependent.

Asymptotic properties of lasso for high-dimensional time series have
been considered by [\citet{powai2012,WuWu2014}], and we provide
detailed comparisons with those studies in Section~\ref{secstoch-reg}.
In short, these works either assume RE conditions or establish their
validity within a very restricted class of $\operatorname{VAR}(1)$ models, as
illustrated in Figure~1
and Lemma E.2
in Appendix~E (supplementary material [\citet{supp}]).

A major contribution of the present study is to establish the validity
of suitable RE and deviation conditions for a large class of stationary
Gaussian processes $\{X^t\}$ and $\{\varepsilon^t\}$. As a result, this
work extends existing results to a much larger class of time series
models and provides deeper insights into the effect of dependence on
the estimation error of lasso.

\textit{Vector autoregression} (VAR) represents a popular class of
time series models in applied macroeconomics and finance, widely used
for structural analysis and simultaneous forecasting of a number of
temporally observed variables [\citet{sims1980,bernanke2005,stock2005FAVAR}]. Unlike structural models, VAR provides a broad
framework for capturing complex temporal and cross-sectional
interrelationship among the time series [\citet{banburra09BVAR}].
In addition to economics, VAR models have been instrumental in linear
system identification problems in control theory [\citet
{kumar1986stochastic}], while more recently, they have become standard
tools in functional genomics for reconstruction of regulatory networks
[\citet{alitrunc,Michailidis-review}] and in neuroscience for
understanding effective connectivity patterns between brain regions
[\citet{Smith-12-neuroimage,friston09,Seth2013}].

Formally, for a $p$-dimensional vector-valued stationary time series $\{
X^t\} = \{(X^t_1, \ldots, X^t_p) \}$, a VAR model of lag $d$ [$\operatorname{VAR}(d)$]
with serially uncorrelated Gaussian errors takes the form
%
\begin{equation}
\label{eqnVAR-defn} X^t = A_1 X^{t-1} + \cdots+
A_d X^{t-d} + \varepsilon^t,\qquad \varepsilon^t \stackrel{\mathrm{i.i.d.}}{\sim} N(\mathbf{0}, \Sigma_\varepsilon),
\end{equation}
where $A_1, \ldots, A_d$ are $p \times p$ matrices and $\varepsilon^t$
is a $p$-dimensional vector of possibly correlated innovation shocks.
The main objective in VAR models is to estimate the transition matrices
$A_1, \ldots, A_d$, together with the order of the model $d$, based on
realizations $\{X^0, X^1, \ldots, X^T \}$. The structure of the
transition matrices provides insight into the complex temporal
relationships amongst the $p$ time series and leads to efficient
forecasting strategies.

VAR estimation is a natural high-dimensional problem, since the
dimensionality of the parameter space ($dp^2$) grows quadratically with
$p$. For example, estimating a $\operatorname{VAR}(2)$ model with $p = 20$ time series
requires estimating $dp^2 = 800$ parameters. However, a comparable
number of stationary observations is rarely available in practice. In
the low-dimensional setting, VAR estimation is carried out by
reformulating it as a multivariate regression problem [\citet
{lutkepohl2005new}]. Under high-dimensional scaling and sparsity
assumptions on the transition matrices, a natural strategy is to resort
to $\ell_1$-penalized least squares or log-likelihood based methods
[\citet{songbickel2011,davis2012}].

Compared to stochastic regression, theoretical analysis of large VAR
requires two important considerations. First, since the response
variable is multivariate, the choice of the loss function (least
squares, negative log-likelihood) plays an important role in estimation
and prediction, especially when the multivariate error process has
correlated components. Second, correlation of the error process with
the process of predictors $\operatorname{Cov}(X^{t}, \varepsilon^{t})
\neq0$ makes the theoretical analysis more involved. Existing work on
high-dimensional VAR models requires stringent assumptions on the
dependence structure [\citet{songbickel2011}], or on the
transition matrix [\citet{negwai2011}], which are violated by
many stable VAR models, as discussed in Section~\ref{secVAR}. Our
results show that consistent estimation is possible with $\ell
_1$-penalization for \textit{both} least squares and log-likelihood
based choices of loss functions under high-dimensional scaling for
\textit{any} stable $\operatorname{VAR}(d)$ models. Interestingly, the latter choice of
loss function leads to an $M$-estimation problem that does not fit into
the stochastic regression framework. As in the case of stochastic
regression, we establish the validity of suitable restricted eigenvalue
and deviation conditions using the stability measures introduced in
this work.

{{A central theme of our theoretical results is that the effect of
dependence on the behavior of these regularized estimates can be nicely
captured by the spectral properties of the underlying multivariate
processes. In particular, we show that the estimation error of lasso in
the time series models scales at the same rate as for i.i.d. data,
modulo a ``price'' of dependence, which can be interpreted as a measure
of ``narrowness'' of the underlying spectra. This agrees with a
fundamental phenomenon in the signal processing literature---a flatter
autocorrelation function (slower decay of temporal dependence)
corresponds to a narrower spectrum and vice versa. Moreover, for linear
ARMA models, our spectral approach has an added advantage of
interpretability, since the spectral density of this class allows a
closed form expression in terms of the model parameters.

At the core of our theoretical results are some novel deviation bounds
for dependent data established in Section~\ref{secmain-results}. These
deviation bounds serve two important purposes. First, they help verify
routinely used restricted eigenvalue and deviation conditions used in
the lasso literature for a large class of time series models and help
develop a theory independent of abstract regularity assumptions.
Second, these deviation bounds are general enough to seamlessly
integrate with the existing theory of other regularization mechanisms
and hence extend the available results to time series setting. Examples
include sparse covariance estimation via hard thresholding, nonconvex
penalties like SCAD and MCP for sparse modeling, group lasso for
structured sparsity and nuclear norm minimization for low-rank
modeling, as discussed in Section~\ref{secdiscussion}. It is worth
noting that many of these regularization mechanisms have been applied
on time series data with good empirical performance [\citet
{songbickel2011,fan-lv-08,BL08covthresh}].
}}

\subsection*{Outline of the paper} The remainder of the paper is organized
as follows. In Section~\ref{secmain-results}, we first demonstrate via
simulation how lasso errors scale in low and high-dimensional regimes
for time series data which motivates the proposed stability measure,
discuss relevant spectral properties of stationary processes, introduce
our measures of stability and present the main deviation bounds used in
subsequent analyses. In Section~\ref{secstoch-reg} we derive
nonasymptotic upper bounds on the estimation error of lasso in
stochastic regression with serially correlated errors. Section~\ref
{secVAR} is devoted to the modeling, estimation and theoretical
analysis of sparse VAR models. We examine both least squares and
likelihood based regularized estimation of VAR models and their
consistency properties. In Section~\ref{secextension}, we discuss
extensions of the current framework to other regularized estimation
problems in high-dimensional time series models.
Finally, Section~\ref{secsim} illustrates the performance of lasso
estimates in stochastic regression and VAR estimation through
simulation studies. We delegate many of the technical proofs to the
Appendices in the supplement [\citet{supp}].


\subsection*{Notation} Throughout this paper, $\mathbb{Z}$, $\mathbb{R}$
and $\mathbb{C}$ denote the sets of integers, real numbers and complex
numbers, respectively. We denote the cardinality of a set by $J$ by
$\llvert  J\rrvert  $. For a vector $v \in\mathbb{R}^p$, we denote $\ell_q$ norms by
$\llVert  v\rrVert  _q:=  (\sum_{j=1}^p \llvert  v_j\rrvert  ^q  )^{1/q}$, for $q > 0$.
We use $\llVert  v\rrVert  _0$ to denote $\llvert  \operatorname{supp}(v)\rrvert   = \sum_{i=1}^p \mathbf{1}[v_j
\neq0]$ and $\llVert  v\rrVert  _{\infty}$ to denote $\max_{j} \llvert  v_j\rrvert  $. Unless
mentioned otherwise, we always use $\llVert  \cdot \rrVert  $ to denote $\ell_2$-norm of
a vector $v$. For a matrix $A$, $\rho(A)$, $\llVert  A\rrVert  $ and $\llVert  A\rrVert  _F$ will
denote its spectral radius $\llvert  \Lambda_{\max}(A)\rrvert  $, operator norm
$\sqrt{\Lambda_{\max}(A'A)}$ and Frobenius norm $\sqrt{\operatorname{tr}(A'A)}$,
respectively. We will also use $\llVert  A\rrVert  _{\max}$, $\llVert  A\rrVert  _1$ and $\llVert  A\rrVert
_{\infty}$ to denote the coordinate-wise maximum (in absolute value),
maximum absolute row sum and maximum absolute column sum of a matrix,
respectively. For any $p \ge1$, $q \ge0$, $r > 0$, we denote the unit
balls by $\mathbb{B}_q(r):= \{v \in\mathbb{R}^p\dvtx  \llVert  v\rrVert  _q \le r \}$.
For any $J \subset\{1, \ldots, p\}$ and $\kappa> 0$, we define the
cone set $\EuScript{C}(S, \kappa) = \{v \in\mathbb{R}^p\dvtx  \llVert  v_{S^c}\rrVert
_1 \le\kappa\llVert  v_S\rrVert  _1 \}$ and the sparse set $\EuScript{K}(s) =
\mathbb{B}_0(s) \cap\mathbb{B}_2(1)$, for any $s \ge1$. For any set
$V$, we denote its closure and convex hull by $\operatorname{cl}\{V\}$ and $\operatorname{conv}\{V\}
$. For a symmetric or Hermitian matrix $A$, we denote its maximum and
minimum eigenvalues by $\Lambda_{\min}(A)$ and $\Lambda_{\max}(A)$.
We use $e_i$ to denote the $i$th unit vector in $\mathbb{R}^p$.
Throughout the paper, we write $A \succsim B$ if there exists an
absolute constant $c$, independent of the model parameters, such that
$A \ge cB$. We use $A \asymp B$ to denote $A \succsim B$ and $B
\succsim A$.

\section{Deviation bounds for multivariate Gaussian time series}\label{secmain-results}

\subsection{Effect of temporal dependence on lasso errors}\label{secmain-results-examples}

Whereas in classical asymptotic analysis of time series, the
quantification of temporal dependence and its impact on the limiting
behavior of the model parameter estimates are typically achieved by
assuming some mixing condition on the underlying stochastic process,
this route is hard to follow in a
high-dimensional context, even for standard ARMA processes. In recent
work, \citet{WuWu2014} and \citet{chen-wu-13-AOS} investigate the
asymptotic properties of lasso and covariance thresholding in the time
series context, assuming a specific rate of decay on the functional
dependence measure [\citet{Wu-PNAS-05}] of the underlying
stationary process. For $\operatorname{VAR}(1)$ processes $X^t=A_1 X^{t-1}+\varepsilon^t$,
the mixing rates and the functional dependence measure are known to
scale with the spectral radius $\rho(A)$ [\citet{liebscher2005,chen-wu-13-AOS}]. The following two simulation experiments show that
dependence in the data affect the convergence rates of lasso estimates
in a more \textit{intricate} manner, not completely captured by $\rho
(A)$. Further, several authors [\citet{powai2012,negwai2011,hanliu13VAR}] conducted nonasymptotic analysis of high-dimensional
$\operatorname{VAR}(1)$ models, assuming $\llVert  A\rrVert   < 1$. In Appendix E (supplementary material [\citet{supp}])
(see Figure~1
and Lemma E.2),
we show that this assumption is restrictive and is violated by many
stable $\operatorname{VAR}(1)$ models. More importantly, such an assumption does not
generalize beyond $\operatorname{VAR}(1)$.

\begin{exa}\label{ex1}
We generate data from the stochastic regression model
(\ref{eqnstoch-reg-defn-intro}) with $p=200$ predictors and i.i.d.
errors $\{\varepsilon^t\}$. The process of predictors comes from a
Gaussian $\operatorname{VAR}(1)$ model $X^t=A X^{t-1}+\xi^t$, where $A$ is an upper
triangular matrix with $\alpha=0.2$ on the diagonal and $\gamma$ on
the two upper off-diagonal bands. We generate processes with different
levels of cross-correlation among the predictors by changing $\gamma$
and plot the average estimation error of lasso (over multiple iterates)
against different sample sizes $n$ in Figure~\ref{figexample-dep}.

%
\begin{figure}

\includegraphics{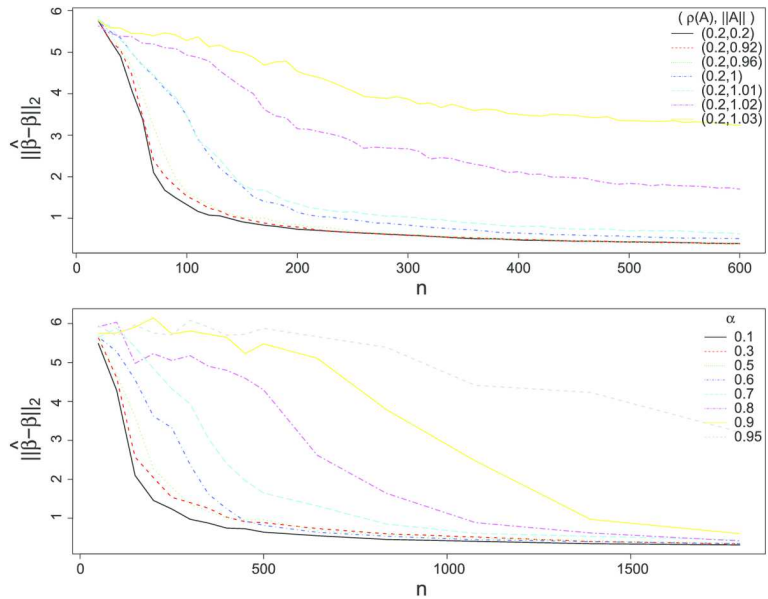}

\caption{Estimation error of lasso in stochastic regression. Top
panel: Example~\protect\ref{ex1}, $\operatorname{VAR}(1)$ process of predictors with cross-sectional
dependence. Bottom panel: Example~\protect\ref{ex2}, $\operatorname{VAR}(2)$ process of predictors with
no cross-sectional dependence.}\label{figexample-dep}
\end{figure}

The spectral radius is \textit{common} ($\alpha=0.2$) across all models.
Consistently with the classical low-dimensional asymptotics, the lasso
errors for different processes seem to converge as $n$ goes to
infinity. However, for small to moderate $n$, as is common in
high-dimensional regimes, lasso errors are considerably different for
different processes. Capturing the effect of cross-dependence via $\llVert
A\rrVert   < 1$ has limitations, as discussed above. We also see that the
errors decay even when $\llVert  A\rrVert  $ exceeds $1$. This motivates a new
approach to capture the cross-dependence among the univariate components.
\end{exa}

\begin{exa}\label{ex2}
Even in the absence of cross-dependence, lasso errors
exhibit interesting behavior in different regimes, as we show in the
next example. Here we generate a similar regression model with $p=500$
predictors, each generated independently from a\vspace*{2pt} Gaussian $\operatorname{VAR}(2)$ process
$X^t_j=2 \alpha X^t_{j-1} - \alpha^2 X^t_{j-2}+\xi^t$, $0 < \alpha<
1$, $\Gamma_X(0)=1$.
The assumption $\llVert  A\rrVert   < 1$ is not applicable here. The processes with
different $\alpha$ exhibit different behavior for small to moderate
$n$, as predicted by their mixing rates and the functional dependence
measures, although it seems the effect of this dependence is
significantly reduced when the sample size is large (Figure~\ref
{figexample-dep}).

These examples motivate us to introduce a different measure to
quantify dependence that reconciles the observed behavior of the lasso errors.
\end{exa}

\subsection{Measure of stability}\label{secmain-results-stability}
Consider a $p$-dimensional discrete time, centered,
covariance-stationary process $\{X^t\}_{t \in\mathbb{Z}}$ with
autocovariance function $\Gamma_X(h) = \operatorname{Cov}(X^t,
X^{t+h})$, $t, h \in\mathbb{Z}$. We make the following assumption:

\begin{assumption}\label{assumpspectral-density}
The spectral density function
%
\begin{equation}
f_{X}(\theta) := \frac{1}{2\pi} \sum_{\ell= -\infty}^{\infty}
\Gamma_X(\ell) e^{-i\ell\theta}, \qquad \theta\in[-\pi, \pi]
\end{equation}
exists, and its maximum eigenvalue is bounded a.e. on $[-\pi, \pi]$,
that is,
%
\begin{equation}
\label{defnmeasures-stability} \mathcal{M}(f_X):= \mathop{\operatorname{ess}
\operatorname {sup}}_{\theta\in[-\pi, \pi]} \Lambda_{\max} \bigl(f_X(
\theta) \bigr) < \infty.
\end{equation}
\end{assumption}

We will often write $f$ instead of $f_X$ and $\Gamma$ instead of
$\Gamma_X$, when the underlying process is clear from the context.
Existence of the spectral density is guaranteed if $\sum_{l =
0}^{\infty}\llVert   \Gamma(l)\rrVert  ^2 < \infty$. Further, if $\sum_{l=0}^\infty\llVert  \Gamma(l)\rrVert   < \infty$, then the spectral density is
bounded, continuous and the essential supremum in the definition of
$\mathcal{M}(f_X)$ is actually the maximum. Assumption \ref
{assumpspectral-density} is satisfied by a large class of general
linear processes, including stable, invertible ARMA processes
[\citet{Priestley2}]. Moreover, the spectral density has a closed
form expression for these processes, as shown in the following examples.

\begin{exa*}
An ARMA($d, \ell$) process $\{X^t \}$
%
\begin{eqnarray}
\label{eqnARMA} X^t &=& A_1 X^{t-1} +
A_2 X^{t-2} + \cdots+ A_d X^{t-d}
\nonumber\\[-8pt]\\[-8pt]\nonumber
&&{} + \varepsilon^t - B_1 \varepsilon^{t-1}
- B_2 \varepsilon^{t-2} - \cdots- B_\ell
\varepsilon^{t-\ell}
\end{eqnarray}
is stable and invertible if the matrix valued polynomials $\mathcal
{A}(z):= I_p - \sum_{t=1}^d A_t z^t$ and $\mathcal{B}(z):= I_p - \sum_{t=1}^{\ell} B_t z^t$ satisfy $\det(\mathcal{A}(z)) \neq0$ and
$\det(\mathcal{B}(z)) \neq0$ on the unit circle of the complex plane
$\{z \in\mathbb{C}\dvtx  \llvert  z\rrvert   = 1 \}$.

For a stable, invertible ARMA process, the spectral density takes the form
%
\begin{equation}
\label{eqnspectral-density-ARMA} f_X (\theta) = \frac{1}{2\pi} \bigl(
\mathcal{A}^{-1} \bigl(e^{-i\theta
}\bigr) \bigr) \mathcal{B}
\bigl(e^{-i\theta}\bigr) \Sigma_{\varepsilon} \mathcal{B}^*
\bigl(e^{-i\theta}\bigr) \bigl( \mathcal{A}^{-1}
\bigl(e^{-i\theta
}\bigr) \bigr)^*.
\end{equation}
In Appendix E (supplementary material [\citet{supp}]),
we provide more details on general linear processes and connection with
mixing conditions.
\end{exa*}

%
\begin{figure}

\includegraphics{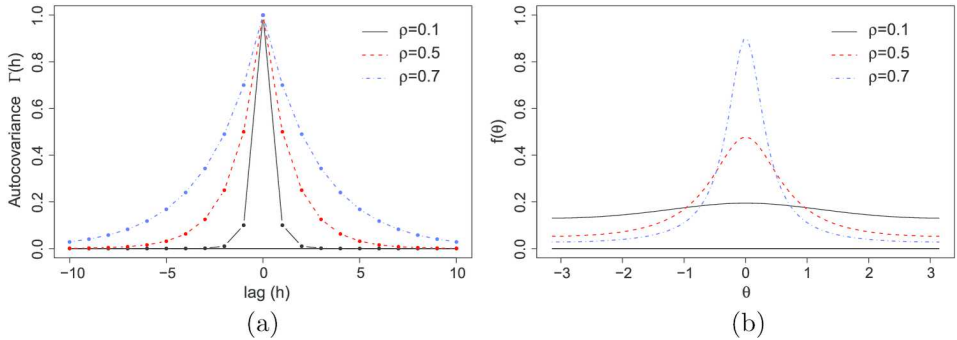}

\caption{Autocovariance $\Gamma(h)$ and spectral density $f(\theta)$
of a univariate $\operatorname{AR}(1)$ process $X^t = \rho X^{t-1} + \varepsilon^t$, $0
< \rho< 1$, $\Gamma_{X}(0) = 1=\int_{-\pi}^{\pi} f(\theta)
\,d\theta$. Processes with stronger temporal dependence, that is, with
larger $\rho$, have flatter $\Gamma$ and narrower $f$. For $\rho=
1$, the process is unstable, and the spectral density does not exist.
\textup{(a)}~Autocovariance of $\operatorname{AR}(1)$,
\textup{(b)}~spectral density of $\operatorname{AR}(1)$.}\label{figGammafAR1}
\end{figure}

Existence of the spectral density ensures the following representation
of the autocovariance matrices
%
\begin{equation}
\label{eqncov-spectral} \Gamma_X(\ell) = \int_{-\pi}^{\pi}
f_X(\theta) e^{i\ell\theta} \,d\theta\qquad\mbox{for all } \ell\in
\mathbb{Z}.
\end{equation}
Since
the autocovariance function characterizes a centered Gaussian process,
it can be used to quantify the temporal and cross-sectional dependence
for this class of models. In particular, spectral density provides
insight into the stability of the process, as illustrated and explained
in the caption of Figure~\ref{figGammafAR1}. The upshot is that the
peak of the spectral density can be used as a measure of stability of
the process.

More generally, for a $p$-dimensional time series $\{X^t\}$, a natural
analogue of the ``peak'' is the maximum eigenvalue of the
(matrix-valued) spectral density function over the unit circle, as
defined in (\ref{defnmeasures-stability}).

In our analysis of high-dimensional time series, we will use $\mathcal
{M}(f_X)$ as a \textit{measure of stability} of the process. Processes
with larger $\mathcal{M}(f_X)$ will be considered less stable.

For any $k$-dimensional subset $J$ of $\{1, \ldots, p\}$, we can
similarly measure the stability of the subprocess $\{X(J)\} = \{
(X^t_j)\dvtx  j \in J\}_{t \in\mathbb{Z}}$ as $\mathcal{M}(f_{X(J)})$. We
will measure the stability of all $k$-dimensional subprocesses of $\{
X^t\}$ using
\[
\label{defnstability-sub-process} \mathcal{M}(f_X, k):= \max_{J \subseteq\{1, \ldots, p \}, \llvert  J\rrvert   \le
k}
\mathcal{M}(f_{X(J)}).
\]
Clearly, $\mathcal{M}(f_X) = \mathcal{M}(f_X, p)$. For completeness,
we define $\mathcal{M}(f_X, k)$ to be $\mathcal{M}(f_X)$, for all $k
\ge p$. It follows from the definitions that
\[
\mathcal{M}(f_X, 1) \le\mathcal{M}(f_X,2) \le\cdots\le
\mathcal {M}(f_X, p) = \mathcal{M}(f_X).
\]
If $\{X^t \}$ and $\{Y^t\}$ are independent $p$-dimensional time series
satisfying Assumption~\ref{assumpspectral-density} and $Z^t =
X^t+Y^t$, then $f_Z = f_X + f_Y$. Consequently,
\[
\label{eqnstability-sum-process} \mathcal{M}(f_Z) = \mathcal{M}(f_X) +
\mathcal{M}(f_Y).
\]
More generally, for any two $p$-dimensional processes $\{X^t\}$ and $\{
Y^t\}$, the cross-spectral density is defined as
\[
f_{X,Y}(\theta)=(1/2\pi) \sum_{l=-\infty}^{\infty}
\Gamma _{X,Y}(l) e^{-il\theta},\qquad \theta\in[-\pi, \pi],
\]
where $\Gamma_{X,Y}(h)=\operatorname{Cov}(X^t, Y^{t+h})$, $h \in
\mathbb{Z}$. If the joint process $W^t=[(X^t)', (Y^t)']'$ satisfies
Assumption \ref{assumpspectral-density}, we can similarly define the
cross-spectral measure of stability
\[
\mathcal{M}(f_{X,Y})=\mathop{\operatorname{ess}\operatorname
{sup}}_{\theta\in[-\pi, \pi]} \sqrt{\Lambda_{\max} \bigl( f^*_{X,Y}(
\theta) f_{X,Y}(\theta) \bigr)}.
\]

For studying stochastic regression and VAR problems, we also need the
lower extremum of the spectral density over the unit circle,
\[
\label{defnmeasures-stability-min} \EuFrak{m}(f_X):=\mathop{\operatorname{ess}
\operatorname {inf}}_{\theta\in[-\pi, \pi]} \Lambda_{\min} \bigl( f_X
(\theta ) \bigr).
\]
Since $\EuFrak{m}(f_X)$ captures the dependence among the univariate
components of the vector-valued time series, it plays a crucial role in
our analysis of high-dimensional regression in quantifying dependence
among the columns of the design matrix.

For stable, invertible ARMA processes and general linear processes with
stable transfer functions, the spectral density is bounded and
continuous. In these cases, the essential supremum (infimum) in the
above definitions of $\EuFrak{m}(f_X)$ and $\mathcal{M}(f_X)$ reduce
to maximum (minimum) because of the continuity of eigenvalues and the
compactness of the unit circle $\{z \in\mathbb{C}\dvtx  \llvert  z\rrvert   = 1 \}$.

Note that $\EuFrak{m}(f_X)$ and $\mathcal{M}(f_X)$ may not have
closed form expressions for general stationary processes. However, for
a stationary ARMA process (\ref{eqnARMA}), we have the following bounds:
%
\begin{eqnarray}
\label{eqnmf-to-mu} \EuFrak{m}(f_X) &\ge&\frac{1}{2\pi}
\frac{\Lambda_{\min}(\Sigma
_{\varepsilon}) \mu_{\min}(\mathcal{B})}{\mu_{\max}(\mathcal{A})},\nonumber
\\
\mathcal{M}(f_X) &\le& \frac{1}{2\pi}
\frac{\Lambda_{\max}(\Sigma
_{\varepsilon}) \mu_{\max}(\mathcal{B})}{\mu_{\min}(\mathcal{A})}
\nonumber\\[-8pt]\\[-8pt]\nonumber
\mu_{\min}(\mathcal{A})&:=& \min_{\llvert  z\rrvert   = 1}
\Lambda_{\min} \bigl(\mathcal {A}^*(z)\mathcal{A}(z)\bigr),
\\
\mu_{\max}(\mathcal{A})&:=& \max_{\llvert  z\rrvert   = 1}
\Lambda_{\max} \bigl(\mathcal {A}^*(z)\mathcal{A}(z)\bigr),
\nonumber
\end{eqnarray}
%
and $\mu_{\min}(\mathcal{B})$, $\mu_{\max}(\mathcal{B})$ are
defined accordingly.

It is often easier to work with $\mu_{\min}(\mathcal{A})$ and $\mu
_{\max}(\mathcal{A})$ instead of $\EuFrak{m}(f_X)$ and $\mathcal
{M}(f_X)$. In particular, we have the following bounds:

\begin{prop}\label{propstability-bounds}
Consider a polynomial $\mathcal{A}(z) = I_p - \sum_{t = 1}^d A_t
z^t$, $z \in\mathbb{C}$, satisfying $\det(\mathcal{A}(z)) \neq0$ for
all $\llvert  z\rrvert   \le1$:
\begin{longlist}[(ii)]
\item[(i)] For any $d \ge1$, $\mu_{\max}(\mathcal{A}) \le [
1 + ({v}_{\mathrm{in}} + {v}_{\mathrm{out}})/2  ]^2$, where
\[
v_{\mathrm{in}} = \sum_{h=1}^d \max
_{1 \le i \le p} \sum_{j=1}^p
\bigl\llvert A_{h}(i,j)\bigr\rrvert,\qquad  v_{\mathrm{out}} = \sum
_{h=1}^d \max_{1 \le j \le p}
\sum_{i=1}^p \bigl\llvert
A_{h}(i,j)\bigr\rrvert .
\]
\item[(ii)] If $d = 1$, and $A_1$ is diagonalizable, then
\[
\mu_{\min}(\mathcal{A}) \ge \bigl(1 - \rho(A_1)
\bigr)^2 \llVert P\rrVert ^{-2} \bigl\llVert
P^{-1}\bigr\rrVert ^{-2},
\]
where $\rho(A_1)$ is the spectral radius (maximum absolute eigenvalue)
of $A_1$, and the columns of $P$ are eigenvectors of $A_1$.
\end{longlist}
\end{prop}

Proposition \ref{propstability-bounds}, together with (\ref
{eqnmf-to-mu}), demonstrate how $\EuFrak{m}(f_X)$ and $\mathcal
{M}(f_X)$ behave for ARMA models. For instance, for a $\operatorname{VAR}(1)$ process,
these quantities are bounded away from zero and infinity as long as the
noise covariance structure and the matrix of eigenvectors of $A_1$ are
well conditioned, the spectral radius of $A_1$ is bounded away from $1$
and the entries of $A_1$ do not concentrate on a single row or column.
The proof is delegated to Appendix E (supplementary material [\citet{supp}]).

\subsection{Deviation bounds} \label{secmain-results-dev-bounds}
Based on realizations of $\{X^t \}_{t=1}^n$ generated according to a
stationary process satisfying Assumption \ref{assumpspectral-density}, we construct the data matrix $\mathcal{X} =
[X^n: \ldots: X^1  ]'$ and the sample Gram matrix $S =
\mathcal{X}'\mathcal{X}/n$. Deriving suitable concentration bounds on
$S$ is a key step for studying regression and VAR estimation problems
in high dimension. In the time series context, this is particularly
challenging, since both the rows and columns of the data matrix
$\mathcal{X}$ are dependent on each other. When the underlying process
is Gaussian, this dependence can be expressed using the covariance
matrix of the random vector $\operatorname{vec}(\mathcal{X}')$. We denote this
covariance matrix by $\Upsilon^{{X}}_n := \operatorname
{Cov}(\operatorname{vec}(\mathcal{X}'), \operatorname{vec}(\mathcal{X}'))_{np \times np}$.

The next proposition provides bounds on the extreme eigenvalues of
$\Upsilon^{{X}}_n$ { and generalizes analogous results in univariate
analysis presented in
\citet{XiaoWu2012} and \citet{grenander1958toeplitz}}. A
similar result {for block Toeplitz forms} under slightly different
conditions can be found in \citet{parter61extemeeigentoeplitz}.
Note that these bounds depend only on the spectral density $f_X$ and
are independent of the sample size $n$.

\begin{prop}\label{proptoeplitz-eigen}
For any $n \ge1$, $p \ge1$,
\[
\label{eqnbounds-eigen} 2\pi\EuFrak{m}(f_X) \le\Lambda_{\min} \bigl(
\Upsilon ^{{X}}_n \bigr) \le\Lambda_{\max} \bigl(
\Upsilon^{{X}}_n \bigr) \le2\pi\mathcal{M}(f_X).
\]
In particular, for $n = 1$,
\[
2\pi\EuFrak{m}(f_X) \le\Lambda_{\min} \bigl(
\Gamma_X(0) \bigr) \le\Lambda_{\max} \bigl(
\Gamma_X(0) \bigr) \le2\pi\mathcal{M}(f_X).
\]
\end{prop}

Next, we establish some deviation bounds on $S=\mathcal{X}'\mathcal
{X}/n$ and $\mathcal{X}'E/n$. These bounds serve as starting points
for analyzing regression and covariance estimation problems. In part
(a), the first deviation bound shows how $\llVert  \mathcal{X}v \rrVert  ^2/ n\llVert  v\rrVert
^2$ concentrates around its expectation, where $v \in\mathbb{R}^p$ is
a fixed vector. This will be used to verify restricted eigenvalue
assumptions for stochastic regression and VAR estimation problems. The
second deviation bound is about the concentration of the entries of $S$
around their expectations. This will be useful for estimating sparse
covariance matrices. In part (b), we establish deviation bounds on how
$\mathcal{X}'\mathcal{Y}/n$ concentrates around zero ($\mathcal{Y}$
is the data matrix from another process $\{Y^t\}$). In regression and
VAR problems, applying this bound with $\{Y^t\}$ as the error process
enables the derivation of necessary deviation bounds on $\mathcal
{X}'E/n$ under different norms.

\begin{prop}\label{propconc-S}
\textup{(a)} For a stationary, centered Gaussian time series $\{X^t\}_{t \in
\mathbb{Z}}$ satisfying Assumption \ref{assumpspectral-density},
there exists a constant $c > 0$ such that for any $k$-sparse vectors
$u, v \in\mathbb{R}^{p}$ with $\llVert  u\rrVert   \le1$, $\llVert  v\rrVert   \le1$, $k \ge
1$, and any $\eta\ge0$,
%
\begin{eqnarray}
\mathbb{P} \bigl[\bigl\llvert v' \bigl( S -
\Gamma_X(0) \bigr) v \bigr\rrvert > 2\pi\mathcal{M}(f_X,
k) \eta \bigr] &\le& 2\exp \bigl[-cn \min\bigl\{ \eta^2, \eta\bigr\}
\bigr], \label{eqnconc-S-op-v}
\\
\mathbb{P} \bigl[\bigl\llvert u' \bigl( S -
\Gamma_X(0) \bigr) v \bigr\rrvert > 6\pi\mathcal{M}(f_X,
2k) \eta \bigr] &\le& 6\exp \bigl[-cn \min\bigl\{ \eta^2, \eta\bigr\}
\bigr]. \label{eqnconc-S-op-uv}
\end{eqnarray}
In particular, for any $i,j \in\{1, \ldots, p \}$, we have
%
\begin{equation}
\label{eqnconc-S-entry} \mathbb{P} \bigl[\bigl\llvert S_{ij} -
\Gamma_{ij}(0) \bigr\rrvert > 6\pi \mathcal{M}(f_X, 2) \eta
\bigr] \le6\exp \bigl[-cn \min\bigl\{ \eta ^2, \eta\bigr\} \bigr].
\end{equation}

\textup{(b)} Consider two $p$-dimensional, centered, stationary Gaussian
processes $\{ X^t\}_{t \in\mathbb{Z}}$ and $\{Y^t\}_{t \in\mathbb
{Z}}$ with $\operatorname{Cov}(X^t, Y^t)=0$ for every $t \in\mathbb
{Z}$ and the joint process $[(X^t)', (Y^t)']'$ satisfying Assumption
\ref{assumpspectral-density}. Let $\mathcal{X}=[X^n: \ldots: X^1]'$
and $\mathcal{Y}=[Y^n:\ldots:Y^1]'$ be the data matrices. Then there
exists a constant $c>0$ such that for any $u,   v \in\mathbb{R}^p$
with $\llVert  u\rrVert  \le1$, $\llVert  v\rrVert  \le1$, we have
%
\begin{eqnarray}
\label{eqnconc-XY}
&& \mathbb{P} \bigl[\bigl\llvert u' \bigl(
\mathcal{X}'\mathcal{Y}/n \bigr) v \bigr\rrvert > 2 \pi \bigl(
\mathcal{M}(f_X)+\mathcal{M}(f_Y)+\mathcal
{M}(f_{X,Y}) \bigr) \eta \bigr]
\nonumber\\[-8pt]\\[-8pt]\nonumber
&&\qquad \le6 \exp \bigl[ -cn \min\bigl\{\eta, \eta^2 \bigr\} \bigr].
\end{eqnarray}
In particular, for any stable $\operatorname{VAR}(d)$ model (\ref{eqnVAR-defn}) with
$\mathcal{X}=[X^n:\ldots:X^1]'$ and $E=[\varepsilon^{n+h}:\ldots
:\varepsilon^{1+h}]'$, $h > 0$, we have
%
\begin{eqnarray}
\label{eqnconc-XY-var}
&& \mathbb{P} \biggl[\bigl\llvert u' \bigl(
\mathcal{X}'E/n\bigr)v \bigr\rrvert > 2\pi \biggl(
\Lambda_{\max}(\Sigma_{\varepsilon}) \biggl(1+\frac{1+\mu_{\max
}(\mathcal{A})}{\mu_{\min}(\mathcal{A})} \biggr)
\biggr) \eta \biggr]
\nonumber\\[-8pt]\\[-8pt]\nonumber
&&\qquad \le6\exp \bigl[-c n \min\bigl\{\eta, \eta^2 \bigr\} \bigr].
\end{eqnarray}
\end{prop}

Next, we give the proofs of the these two key propositions that employ
techniques in spectral theory of multivariate time series and
nonasymptotic random matrix theory results.
\begin{pf*}{Proof of Proposition \ref{proptoeplitz-eigen}}
For $1 \le r, s \le n$, the $(r,s)$th block of the $np \times np$
matrix $\Upsilon^{{X}}_n$ is a $p \times p$ matrix
\[
\Gamma_X(r-s) = \operatorname{Cov} \bigl(X^{n-r+1},
X^{n-s+1} \bigr).
\]
For any $x \in\mathbb{R}^{np}$, $\llVert  x \rrVert   = 1$, write $x$ as $x = \{
(x^1)', (x^2)', \ldots, (x^{p})' \}'$, where each $x^i \in\mathbb
{R}^p$. Define $G(\theta) = \sum_{r=1}^n x^r e^{-ir\theta}$, for
$\theta\in[-\pi, \pi]$. Note that
%
\begin{eqnarray}
\label{eqnnorm-Gtheta} \int_{-\pi}^{\pi} G^*(\theta) G(\theta)
\,d\theta&=& \sum_{r=1}^n \sum
_{s=1}^n \int_{-\pi}^\pi
\bigl(x^r\bigr)'\bigl(x^s\bigr)
e^{i(r-s)\theta} \,d\theta
\nonumber\\[-8pt]\\[-8pt]\nonumber
&=& \sum_{r=1}^n \bigl\llVert
x^r\bigr\rrVert ^2 2\pi= 2\pi.
\end{eqnarray}
Also,
\begin{eqnarray*}
x' \Upsilon^{{X}}_n x &=& \sum
_{r=1}^n \sum_{s=1}^n
\bigl(x^r\bigr)' \Gamma _X(r-s)
\bigl(x^s\bigr)
\\
&=& \sum_{r=1}^n \sum
_{s=1}^n \int_{-\pi}^\pi
\bigl(x^r\bigr)' f_X(\theta)
e^{i(r-s)\theta} \bigl(x^s\bigr) \,d\theta\qquad\mbox{using (\ref{eqncov-spectral})}
\\
&=& \int_{-\pi}^\pi G^*(\theta) f_X(
\theta) G(\theta) \,d\theta.
\end{eqnarray*}
Since $f_X(\theta)$ is Hermitian, $G^*(\theta) f_X(\theta) G(\theta
)$ is real, for all $\theta\in[-\pi, \pi]$, and
\[
\EuFrak{m}(f_X) G^*(\theta) G(\theta) \le G^*(\theta)
f_X(\theta) G(\theta) \le\mathcal{M}(f_X) G^*(\theta)
G(\theta).
\]
This, together with (\ref{eqnnorm-Gtheta}), implies
\[
2\pi\EuFrak{m}(f_X) \le x' \Upsilon^{{X}}_n
x \le2\pi\mathcal{M}(f_X)
\]
for all $x \in\mathbb{R}^{np}$, $\llVert  x \rrVert   = 1$.
\end{pf*}
%
%
\begin{pf*}{Proof of Proposition \ref{propconc-S}}
(a) First, note that it is enough to prove (\ref{eqnconc-S-op-v}) for
$\llVert  v\rrVert   = 1$.
For any $v \in\mathbb{R}^p$, $\llVert  v\rrVert  =1$, let $J$ denote its support
$\operatorname{supp}(v)$ so that $\llvert  J\rrvert   = k$. define $Y = \mathcal{X}v = \mathcal{X}_J
v_J$. Then $Y \sim N(0_{n\times1}, Q_{n \times n})$ with
\[
Q_{rs} = v_J' \operatorname{Cov}
\bigl(X_J^{n-r+1}, X_J^{n-s+1}\bigr)
v_J = v_J' \Gamma_{X(J)}(r-s)
v_J\qquad\mbox{for all $1\le r, s \le n$}.
\]
Note that $v'Sv=(1/n) Y'Y=(1/n)Z'QZ$ where $Z \sim N(0, I_n)$. Also,\break 
$v'\Gamma_X(0)v=v_J' \Gamma_{X(J)}(0)v_J = \mathbb{E}[Z'QZ/n]$.

So, by the Hanson--Wright inequality of \citet
{RV13-hansen-wright}, with $\llVert  Z_i\rrVert  _{\psi_2} \le1$ since $Z_i \sim
N(0,1)$, we get
%
\begin{eqnarray}
\label{eqnconc-S-hw} \mathbb{P} \bigl[ \bigl\llvert v' \bigl(S-
\Gamma_X(0) \bigr)v\bigr\rrvert > \zeta \bigr] &=& \mathbb{P} \bigl[
\bigl\llvert Z'QZ-\mathbb{E}\bigl[Z'QZ\bigr] \bigr
\rrvert > n\zeta \bigr]
\nonumber\\[-8pt]\\[-8pt]\nonumber
&\le& 2 \exp \biggl[ -cn \min \biggl\{ \frac{n^2 \zeta^2}{\llVert  Q\rrVert
_{F}^2}, \frac{n \zeta}{\llVert  Q\rrVert  }
\biggr\} \biggr].
\end{eqnarray}
Since $\llVert  Q\rrVert  _{F}^2/n \le\llVert  Q\rrVert  ^2$, setting $\zeta=\llVert  Q\rrVert   \eta$, we obtain
\[
\mathbb{P} \bigl[\bigl\llvert v'\bigl(S-\Gamma_X(0)
\bigr)v\bigr\rrvert > \eta\llVert Q\rrVert \bigr] \le2 \exp \bigl[-cn \min\bigl\{
\eta, \eta^2 \bigr\} \bigr].
\]
%
Also, for any $w \in\mathbb{R}^n$, $\llVert  w\rrVert   = 1$, we have
\begin{eqnarray*}
w'Q w &=& \sum_{r=1}^{n}
\sum_{s=1}^n w_r
w_s Q_{rs} = \sum_{r=1}^n
\sum_{s=1}^n w_r
w_s v_J'\Gamma_{X(J)}(r-s)v_J
\\
&=& (w \otimes v)' \Upsilon^{{X(J)}}_n (w
\otimes v)
\\
&\le& \Lambda_{\max} \bigl(\Upsilon^{{X(J)}}_n \bigr)\qquad\mbox{since $\llVert w \otimes v \rrVert = 1$}
\\
&\le& 2\pi\mathcal{M}(f_{X(J)})\qquad \mbox{by Proposition \ref
{proptoeplitz-eigen}}
\\
&\le& 2\pi\mathcal{M}(f_X, k).
\end{eqnarray*}
This establishes an upper bound on the operator norm $\llVert  Q\rrVert   \le2\pi
\mathcal{M}(f_X, k)$.

To prove (\ref{eqnconc-S-op-uv}), note that
\begin{eqnarray*}
2 \bigl\llvert u' \bigl(S - \Gamma_X(0) \bigr) v
\bigr\rrvert &\le& \bigl\llvert u' \bigl(S - \Gamma_X(0)
\bigr)u \bigr\rrvert + \bigl\llvert v' \bigl(S -
\Gamma_X(0)\bigr)v \bigr\rrvert
\\
\nonumber
&&{} + \bigl\llvert (u+v)' \bigl(S - \Gamma_X(0)
\bigr) (u+v) \bigr\rrvert
\end{eqnarray*}
and $u+v$ is $2k$-sparse with $\llVert  u +v\rrVert   \le2$. The result follows by
applying (\ref{eqnconc-S-op-v}) separately on each of the three terms
on the right.

The element-wise deviation bound (\ref{eqnconc-S-entry}) is obtained
by choosing $u = e_i$, $v = e_j$.

(b) Note that $u'(\mathcal{X}'\mathcal{Y}/n)v$ can be viewed as
$(1/n)\sum_{t=1}^n w^t z^t$, where $w^t= \langle u, X^t \rangle$,
$z^t= \langle v, Y^t \rangle$ are two univariate stationary processes
with spectral densities $f_w(\theta) = u'f_X(\theta) u$ and
$f_z(\theta) = v' f_Y(\theta) v$. Since $\operatorname{Cov}(w^t,
z^t)=0$, we have the following decomposition:
\begin{eqnarray*}
\frac{2}{n}\sum_{t=1}^n
w^t z^t &=& \Biggl[ \frac{1}{n} \sum
_{t=1}^n \bigl(w^t+z^t
\bigr)^2-\operatorname{Var}\bigl(w^1+ z^1
\bigr) \Biggr]
\\
&&{} - \Biggl[ \frac{1}{n} \sum_{t=1}^n
\bigl(w^t\bigr)^2-\operatorname {Var}
\bigl(w^1\bigr) \Biggr] - \Biggl[ \frac{1}{n} \sum
_{t=1}^n \bigl(z^t\bigr)^2-
\operatorname{Var}\bigl(z^1\bigr) \Biggr] ,
\end{eqnarray*}
and it suffices to concentrate the three terms separately. Applying
(\ref{eqnconc-S-op-v}) on the process $w^t= \langle u, X^t \rangle$
and noting that $\mathcal{M}(f_w) \le\mathcal{M}(f_X)$, we have
\begin{eqnarray*}
&& \mathbb{P} \Biggl[ \Biggl\llvert (1/n)\sum_{t=1}^n
\bigl(w^t\bigr)^2-\operatorname {Var}
\bigl(w^1\bigr) \Biggr\rrvert >2\pi\mathcal{M}(f_X) \eta
\Biggr]> 2 \exp \bigl[-cn \min\bigl\{\eta, \eta^2 \bigr\} \bigr].
\end{eqnarray*}
A similar argument for $\{z^t \}$ leads to
\begin{eqnarray*}
&& \mathbb{P} \Biggl[ \Biggl\llvert (1/n)\sum_{t=1}^n
\bigl(z^t\bigr)^2-\operatorname {Var}
\bigl(z^1\bigr) \Biggr\rrvert >2\pi\mathcal{M}(f_Y) \eta
\Biggr]> 2 \exp \bigl[-cn \min\bigl\{\eta, \eta^2 \bigr\} \bigr].
\end{eqnarray*}
To concentrate the first term, note that the process $\{w^t+z^t\}$ has
a spectral density given by
\begin{eqnarray*}
f_{w+z}(\theta) &=&\lleft[ \matrix{ u' &
v'} \rright] \lleft[ \matrix{ f_X(\theta) &
f_{X, Y}(\theta)
\vspace*{3pt}\cr
f^*_{X, Y}(\theta) & f_{Y}(
\theta)} \rright] \lleft[ \matrix{ u
\cr
v} \rright]
\\
&=& u'f_X(\theta) u + v'
f_Y(\theta)v + u'f_{X,Y}(\theta) v +
v' f_{X, Y}^*(\theta) u.
\nonumber
\end{eqnarray*}
Since $\llVert  u\rrVert  \le1$, $\llVert  v\rrVert  \le1$, $\mathcal{M}(f_{w+z}) \le\mathcal
{M}(f_X)+\mathcal{M}(f_Y)+2\mathcal{M}(f_{X, Y})$, where the last
term is obtained by applying the Cauchy--Schwarz inequality on each of
the cross-product terms. Applying (\ref{eqnconc-S-op-v}) separately on
$\{w^t\}$, $\{z^t\}$ and $\{w^t+z^t\}$ with the above bounds on the
respective stability measures leads to the final result.

In the special case of a $\operatorname{VAR}(d)$ process, set $\tilde{\varepsilon}^t :=
\varepsilon^{t+h}$ so that\break \mbox{$\operatorname{Cov}(X^t, \tilde{\varepsilon
}^t) = 0$}. Then it suffices to establish upper bounds on $\mathcal
{M}(f_X)$, $\mathcal{M}(f_{\tilde{\varepsilon}})$ and $\mathcal
{M}(f_{X, \tilde{\varepsilon}})$. From~(\ref{eqnmf-to-mu}), $2\pi
\mathcal{M}(f_X)$ is upper bounded by $\Lambda_{\max}(\Sigma
_{\varepsilon})/\mu_{\min}(\mathcal{A})$. The process $\{\tilde
{\varepsilon}^t\}$ is serially uncorrelated, so $\mathcal{M}(f_{\tilde
{\varepsilon}})$ is the same as $\Lambda_{\max}(\Sigma_{\varepsilon})$.
To derive an upper bound on the cross-spectral measure of stability,
note that
\begin{eqnarray*}
\operatorname{Cov}\bigl(X^t, \varepsilon^{t+h+l}\bigr) &=&
\operatorname {Cov} \bigl( X^t, X^{t+h+l}-A_1
X^{t+h+l-1}-\cdots- A_d X^{t+h+l-d} \bigr)
\\
&=& \Gamma_X(h+l)-\Gamma_X(h+l-1) A_1'-
\cdots- \Gamma_X(h+l-d) A_d'.
\end{eqnarray*}
Hence, the cross-spectrum of $\{X^t\}$ and $\{\tilde{\varepsilon}^{t}\}$
can be expressed as
\begin{eqnarray*}
&& f_{X, \tilde{\varepsilon}}(\theta)
\\
&&\qquad = \frac{1}{2\pi} \sum
_{l=-\infty}^{\infty} \bigl[ \Gamma_X(h+l)-
\Gamma_X(h+l-1) A_1'-\cdots-
\Gamma_X(h+l-d) A_d' \bigr]
e^{-il\theta}
\\
&&\qquad = f_X(\theta)e^{ih\theta} \bigl[I-A_1'
e^{-i\theta}-\cdots- A_d' e^{-id\theta} \bigr]
\\
&&\qquad = e^{ih\theta} f_X(\theta) \mathcal{A}^*\bigl(e^{i\theta}
\bigr).
\end{eqnarray*}
Hence $\mathcal{M}(f_{X,\tilde{\varepsilon}})$ is bounded above by
$\mathcal{M}(f_X) \mu_{\max}(\mathcal{A})$. Combining the three
upper bounds on the stability measures and replacing $\mathcal
{M}(f_X)$ with its upper bound in~(\ref{eqnmf-to-mu}), establishes the
final result.
\end{pf*}

\subsubsection*{Role of the two tails in (\protect\ref{eqnconc-S-hw}) and sharpness of the bounds}
The convergence rates of lasso and other regularized estimates in
high-dimensional settings depend on how $S$ concentrates around $\Gamma
_X(0)$ and $\mathcal{X}'E/n$ around $0$, as is evident in subsequent
proofs. In the bounds established above, the effect of dependence is
captured by $\mathcal{M}(f_X)$. In the special case of no temporal and
cross-sectional dependence, our results recover the bounds of lasso for
i.i.d. data, as we remark in Section~\ref{secstoch-reg}. For processes
with strong dependence, however, we believe this bound can be further
sharpened, although a closed form solution of the exact rate was not
established. Next, we provide an asymptotic argument for a fixed $p$
case and demonstrate that in a low-dimensional setting with very large
sample sizes, the effect of dependence can be captured by the
integrated spectrum, which provides a tighter bound.

The sub-Gaussian and sub-exponential tails in the main concentration
inequality (\ref{eqnconc-S-hw}) suggest an interesting phenomenon,
that temporal dependence in the data may affect the concentration
property and in turn the convergence rates of the regularized estimates
in two different ways, depending on which term in the tail bound is dominant.

In the special case of no temporal dependence, that is, $X^t \stackrel
{\mathrm{i.i.d.}}{\sim} N(0, \Sigma)$, the matrix $Q$ is diagonal and $\llVert  Q\rrVert
_{F}/\sqrt{n} = \llVert  Q\rrVert  $. So, setting $\zeta=\eta\llVert  Q\rrVert  _{F} / \sqrt
{n}$ or $\zeta=\eta\llVert  Q\rrVert  $ leads to the same bound, and we recover the
Bernstein-type tail bounds for subexponential random variables
[\citet{vershynin2010-rmt}].

\begin{figure}

\includegraphics{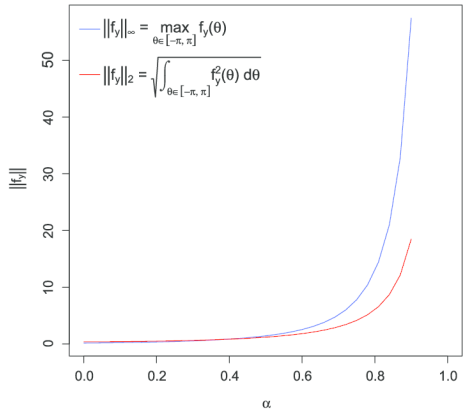}

\caption{$\llVert  f_y \rrVert  _2$ and $\llVert  f_y\rrVert  _\infty$ for a univariate Gaussian
$\operatorname{AR}(2)$ process $y^t = 2\alpha y^{t-1} - \alpha^2 y^{t-2}+\xi^t$,
$\Gamma_y(0)=1$, $0 < \alpha< 1$.}
\label{figspectrum-sharp-picture}
\end{figure}

In the presence of temporal dependence, the two norms $\llVert  Q\rrVert  _{F}$ and
$\llVert  Q\rrVert  $ behave differently, and this affects the rates. To illustrate
this further, we need additional notation. First note that $\mathcal
{M}(f_X)$ can be viewed as $\sup_{\llVert  v\rrVert  =1} \llVert   f_y\rrVert
_\infty$ where $y^t = \langle v, X^t \rangle$ and $\llVert  \cdot \rrVert  _{\infty}$
denotes the $L_\infty$ or sup norm of a function. A related quantity
that will be useful for studying the tails is the Euclidean or $L_2$
norm $\llVert  f_y\rrVert  _2= (\int_{-\pi}^{\pi} f_y^2(\theta) \,d\theta
 )^{1/2}$. For any univariate Gaussian process $\{y^t\}$, it is
easy to see that $\llVert  f_y\rrVert  _2 \le\sqrt{2\pi} \llVert  f_y\rrVert  _{\infty}$, and
they coincide when the process is serially uncorrelated, that is, the
spectrum is flat a.e. With stronger temporal dependence, the spectrum
becomes more spiky and $\llVert  f_y\rrVert  _{\infty}$ changes more sharply\vspace*{2pt} than $\llVert
f_y\rrVert  _2$. In Figure~\ref{figspectrum-sharp-picture}, we demonstrate
this on a family of $\operatorname{AR}(2)$ processes $y^t = 2\alpha y^{t-1} - \alpha^2
y^{t-2}+\xi^t$, $\Gamma_y(0)=1$, $0 < \alpha< 1$.

%

Coming back to the behavior of the two tails, note that
\[
\mathbb{P} \bigl[ \bigl\llvert v'\bigl(S-\Gamma(0)\bigr)v\bigr
\rrvert > \zeta \bigr] \le2 \exp \biggl[-c \min \biggl\{\frac{n \zeta^2}{\llVert  Q\rrVert  _{F}^2/n},
\frac{n\zeta
}{\llVert  Q\rrVert  } \biggr\} \biggr].
\]

We consider a low-dimensional, fixed $p$ regime. It is known that [cf.
Chapter $5$, \citet{grenander1958toeplitz}] for large\vspace*{1pt} $n$, $\llVert  Q\rrVert
^2_{F}/n$ approaches $2\pi\llVert  f_y\rrVert  _2^2$ and $\llVert  Q\rrVert  $ approaches $2\pi\llVert
f_y\rrVert  _{\infty}$. With a choice of $\zeta\asymp\sqrt{\log p /n}$,
the tail probability on the right-hand side can be approximated by
\[
2 \exp \biggl[-c \min \biggl\{\frac{\log p}{c_1 \llVert  f_y \rrVert  _{2}^2}, \frac
{\sqrt{n \log p}}{\llVert  f_y\rrVert  _{\infty}} \biggr\}
\biggr].
\]
This indicates that for very large $n$, the first term will be smaller,
and the tail probability will scale with $\llVert  f_y\rrVert  _2$. So processes with
various levels of dependence should behave similarly in terms of
estimation errors. For strongly dependent processes, where $\Vert f_y\Vert _2
\ll\Vert f_y\Vert _{\infty}$, it would take more samples $n$ for the first
term to offset the second term. With a smaller sample size, the tail
behavior will be driven by $\llVert  f_y\rrVert  _{\infty}$, and the effect of
dependence will be more prominent in the estimation error of the
regularized estimates. Interestingly, this is the same pattern
reflected in Figure~\ref{figexample-dep}.

\section{Stochastic regression}\label{secstoch-reg}

In the presence of serially correlated errors, and under a sparsity
assumption on $\beta^*$, we use the deviation bounds of Section~\ref
{secmain-results} to derive an upper bound on the estimation error of
lasso. Our results show that consistent estimation of $\beta^*$ is
possible, as long as the predictor and noise processes are stable.
We consider the lasso estimate (\ref{eqnlasso-intro}) for the
stochastic regression model (\ref{eqnstoch-reg-defn-intro}). Further,
we assume that both $f_X$ and $f_\varepsilon$ satisfy Assumption \ref
{assumpspectral-density}, and $\beta^*$ is $k$-sparse, with support
$J$, that is, $\llvert  J\rrvert   = k$.

Note that in the low-dimensional regime, consistent estimation relies
on the following assumptions:
\begin{longlist}[(a)]
\item[(a)] $\mathcal{X}' \mathcal{X}/n$ converges to a nonsingular
matrix ($\lim_{n \rightarrow\infty} \Lambda_{\min}  (\frac
{\mathcal{X}'\mathcal{X}}{N}  ) > 0$).
\item[(b)] $\mathcal{X}'E/n$ converges to zero.
\end{longlist}
In the high-dimensional regime ($n \ll p$), the first assumption is
never true since the design matrix is rank-deficient (i.e., more
variables than observations). The second assumption is also very
stringent, since the dimension of $\mathcal{X}'E$ grows with $n$ and $p$.
Interestingly, consistent estimation in the high-dimensional regime can
be ensured under two analogous sufficient conditions. The first one
comes from a class of conditions commonly referred to as \textit{restricted eigenvalue} (RE) conditions [\citet
{bickel2009simultaneous,vandegeerconditions}].\vspace*{1pt} Roughly speaking, these
assumptions require that $\llVert  \mathcal{X}(\hat{\beta} - \beta^*)\rrVert  $
is small only when $\llVert  \hat{\beta} - \beta^*\rrVert  $ is small. For sparse
$\beta^*$ and $\lambda_n$ appropriately chosen, it is now well
understood that the vectors $v = \hat{\beta} - \beta^*$ only vary on
a small subset of the high-dimensional space $\mathbb{R}^p$
[\citet{negahban2012unified}]. As shown in the proof of
Proposition \ref{propmain-result-reg}, the error vectors $v$ in
stochastic regression lie in a cone
\[
\EuScript{C}(J, 3) = \bigl\{ v \in\mathbb{R}^p\dvtx  \llVert
v_{J^c} \rrVert _1 \le 3 \llVert v_J \rrVert
_1 \bigr\},
\]
whenever $\lambda_n \ge4\llVert  \mathcal{X}'E/n \rrVert  _\infty$. This
indicates that the RE condition may not be very stringent after all,
even though $\mathcal{X}$ is singular. Note that verifying that the
assumption indeed holds with high probability is a nontrivial task.

The next proposition shows that a restricted eigenvalue (RE) condition
holds with high probability when the sample size is sufficiently large
and the process of predictors $\{X^t \}$ is stable, with a full-rank
spectral density.

\begin{prop}[(Restricted eigenvalue)]\label{propRE-reg}
If $\EuFrak{m}(f_X) > 0$, then there exist constants $c_i>0$ such that for
$n \succsim\max\{1, \omega^2\}   \min\{k \log(c_0 p/k), k \log p\}$,
\begin{eqnarray*}
&&\mathbb{P} \biggl[ \inf_{v \in\EuScript{C}(J, 3)\backslash\{0\}} \frac{\llVert  \mathcal{X}v\rrVert  ^2}{n \llVert  v\rrVert  ^2} \ge
\alpha_{\mathrm{RE}} \biggr] \ge1 - c_1 \exp \bigl[-c_2
n \min\bigl\{1, \omega^{-2}\bigr\} \bigr],
\end{eqnarray*}
where $\alpha_{\mathrm{RE}} = \pi\EuFrak{m}(f_X)$, $\omega= c_3 \mathcal
{M}(f_X, 2k)/\EuFrak{m}(f_X)$.
\end{prop}

\begin{rems*}
(a) The assumption $\EuFrak{m}(f_X) > 0$ is fairly mild and holds for
stable, invertible ARMA processes. However, the conclusion holds under
weaker assumptions like $\Lambda_{\min}(\Gamma_X(0)) > 0$ or an RE
condition on $\Gamma_X(0)$, replacing $2\pi\EuFrak{m}(f_X)$ by the
minimum (or restricted) eigenvalue of $\Gamma_X(0)$, as evident in the
proof of this proposition.

(b) For large $k$, $k \log(c_0p/k)$ can be much smaller than $k \log
p$, the sample size required for consistent estimation with lasso.

(c) The factor $\omega\asymp\mathcal{M}(f_X, 2k)/\EuFrak{m}(f_X)$
captures the effect of temporal and cross-sectional dependence in the
data. Larger values of $\mathcal{M}(\cdot)$ and smaller values of $\EuFrak
{m}(\cdot)$ indicate stronger dependence in the data, and { the bounds
indicate that} more samples are required to ensure RE holds with high
probability. We demonstrate this on three special types of dependence
in the design matrix $\mathcal{X}$, independent entries, independent
rows and independent columns:
\begin{longlist}[(iii)]
\item[(i)] If the entries of $\mathcal{X}$ are independent from a
$N(0, \sigma^2)$ distribution, we have $\Gamma_X(0) = \sigma^2I$ and
$\Gamma_X(h) = \mathbf{0}$ for $h \neq0$. In this case, $f_X(\theta
) \equiv(1/2\pi)  \sigma^2 I$ and $\mathcal{M}(f_X, 2k)/\EuFrak
{m}(f_X) = 1$.
\item[(ii)] If the rows of $\mathcal{X}$ are independent and
identically distributed as $N(0,\break  \Sigma_{X})$, that is, $\Gamma_X(0)
= \Sigma_{X}$, $\Gamma_X(h) = \mathbf{0}$ for $h \neq0$, the
spectral density takes the form $f_X(\theta) \equiv(1/2\pi) \Sigma
_X$, and $\mathcal{M}(f_X, 2k)/\EuFrak{m}(f_X)$ can be at most
$\Lambda_{\max} (\Sigma_X) /\break  \Lambda_{\min} (\Sigma_X)$.
\item[(iii)] If the columns of $\mathcal{X}$ are independent, that
is, all the univariate components of $\{X^t\}$ are independently
generated according to a common stationary process with spectral
density $f$, then the spectral density of $\{X^t\}$ is $f_X(\theta) =
f(\theta)  I$, and we have
\[
\mathcal{M}(f_X, 2k)/\EuFrak{m}(f_X) = \max
_{\theta\in[-\pi, \pi
]} f(\theta) / \min_{\theta\in[-\pi, \pi]} f(\theta).
\]
The ratio on the right can be viewed as a measure of narrowness of $f$.
Since narrower spectral densities correspond to processes with flatter
autocovariance, this indicates that more samples are needed when the
dependence is stronger.
\end{longlist}

The second sufficient condition for consistency of lasso requires that
the coordinates of $\mathcal{X}'E/n$ uniformly concentrate around $0$.
In the next proposition, we establish a deviation bound on $\llVert
\mathcal{X}'E/n \rrVert  _{\infty}$ that holds with high probability.
Similar results were established in \citet{powai2012} for a $\operatorname{VAR}(1)$
process with serially uncorrelated errors, under the assumption $\llVert  A_1\rrVert
<1$. Our result relies on different techniques, holds for a much larger
class of stationary processes and allows for serial correlation in the
noise term, as well.
\end{rems*}

%
\begin{prop}[(Deviation condition)]\label{propdevn-reg}
For $n \succsim\log p$, there exist constants $c_i > 0$ such that
\[
\mathbb{P} \biggl[\frac{1}{n}\bigl\llVert \mathcal{X}'E
\bigr\rrVert _{\infty
} > c_0 2\pi \bigl[\mathcal{M}(f_X,
1) + \mathcal{M}(f_\varepsilon) \bigr]\sqrt{\frac{\log p}{n}} \biggr] \le
c_1 \exp [ -c_2 \log p ] .
\]
\end{prop}

\begin{rem*}
The deviation inequality shows that the coordinates
of $\mathcal{X}'E/n$ uniformly concentrate around $0$, as long as
the stability measures of $\{\varepsilon^t\}$ and the univariate
components of $\{X^t\}$ grow at a rate slower than $\sqrt{n/\log p}$.
These two propositions allow us to establish error rates for estimation
and prediction in stochastic regression.
\end{rem*}

%
\begin{prop}[(Estimation and prediction error)]\label{propmain-result-reg}
Consider the\break stochastic regression setup of (\ref
{eqnstoch-reg-defn-intro}). If $\beta^*$ is $k$-sparse, $n \succsim
[\mathcal{M}(f_X, k)/\break \EuFrak{m}(f_X)]^2 k \log p$, then there exist
constants $c_i > 0$ such that for
\[
\lambda_n \ge c_0 2\pi \bigl[\mathcal{M}(f_X,
1) + \mathcal {M}(f_\varepsilon) \bigr] \sqrt{(\log p) / n},
\]
any solution $\hat{\beta}$ of (\ref{eqnlasso-intro}) satisfies, with
probability at least $1 - c_1 \exp [-c_2\log p ]$,
\begin{eqnarray*}
\bigl\llVert \hat{\beta} - \beta^* \bigr\rrVert &\le& \frac{2\lambda
_n\sqrt{k}}{\alpha_{\mathrm{RE}}},
\\
\bigl\llVert \hat{\beta} - \beta^* \bigr\rrVert _1 &\le&
\frac{8\lambda_n
k}{\alpha_{\mathrm{RE}}},
\\
\frac{1}{n} \bigl\llVert \mathcal{X} \bigl(\hat{\beta} - \beta^*\bigr)
\bigr\rrVert ^2 &\le& \frac{4 \lambda_n^2 k}{\alpha_{\mathrm{RE}}},
\end{eqnarray*}
where the restricted eigenvalue $\alpha_{\mathrm{RE}} = \pi\EuFrak{m}(f_X)$.

Further, a thresholded variant of lasso $\tilde{\beta}$, defined as
$\tilde{\beta}_j = \{\hat{\beta}_j \mathbf{1}_{\llvert  \hat{\beta
}_j\rrvert   > \lambda_n} \}$, for $1 \le j \le p$, satisfies, with the
same probability,
%
\begin{equation}
\label{eqnthresholded-lasso} \bigl\llvert \operatorname{supp}(\tilde{\beta}) \backslash \operatorname{supp}\bigl(\beta^*
\bigr) \bigr\rrvert \le \frac{24k}{\alpha_{\mathrm{RE}}}.
\end{equation}
\end{prop}

\begin{rems*}
(a) The convergence rates of $\ell_2$-estimation and prediction $\sqrt
{k \log p   /n}$ are of the same order as the rates for regression
with i.i.d. samples. The temporal dependence contributes the additional
term $ [\mathcal{M}(f_X, 1) + \mathcal{M}(f_\varepsilon)
]/\EuFrak{m}(f_X)$ in the error rates and $[\mathcal{M}(f_X,
2k)/\EuFrak{m}(f_X)]^2$ in the sample size requirement. This ensures
fast convergence rates of lasso under high-dimensional scaling, as long
as the processes of predictors and noise are stable.

(b) A thresholded version of lasso enjoys small false positive rates,
as shown in~(\ref{eqnthresholded-lasso}). Note that we do not assume
any ``beta-min'' condition, that is, a lower bound on the minimum signal
strength. It is possible to control the false negatives under suitable
``beta-min'' conditions, as shown in [\citet{zhou2010thresholded}].
\end{rems*}

\subsection*{Comparison with existing results} The problem of stochastic
regression in a high-dimensional setting has been addressed by
\citet{powai2012}. After initial submission of this work, we
became aware of a recent work by \citet{WuWu2014}. Next, we
briefly illustrate the major differences of our results with these
other studies. \citet{powai2012} assume that the process of
predictors $\{X^t \}$ follows a Gaussian $\operatorname{VAR}(1)$ process with transition
matrix satisfying $\llVert  A\rrVert  <1$. They also assume that the errors are
independent. Our results allow both the predictors and the errors to be
generated from any stable Gaussian process. \citet{WuWu2014}
consider lasso estimation with a fixed design matrix and assume that an
RE condition is satisfied. In our work, we consider a random Gaussian
design and establish that RE holds with high probability for a large
class of stable processes. Consequently, our final results of
consistency do not rely on any RE type assumptions. \citet
{WuWu2014} also consider random design regression using a CLIME
estimator and provide an upper bound on the estimation error, without
assuming RE type conditions. However, the established upper bounds seem
to worsen with stronger signal ($\llvert  \beta\rrvert  _1$). Our results do not
exhibit such properties. Finally, both these papers consider a
short-range dependence regime, although their results are derived under
a mild moment condition on the random variables while we focus on
Gaussian processes only. The results in the above paper quantify
dependence via the functional and predictive measure of \citet
{Wu-PNAS-05} and assume a certain decay condition on this measure. For
the multivariate stationary linear processes, this is verified under
another decay condition on the transition matrices in its AR
representation [\citet{chen-wu-13-AOS}]. Our results, on the
other hand, rely on existence and boundedness of spectral density, and
this assumption is satisfied by commonly used stable processes,
including ARMA and general linear processes.

\section{Transition matrix estimation in sparse VAR models}\label{secVAR}

This problem has been considered by several authors in recent years
[\citet{songbickel2011,davis2012,hanliu13VAR}]. Most of these
studies consider a least squares based objective function or estimating
equation to obtain the estimates, which is agnostic to the presence of
cross-correlations among the error components (nondiagonal $\Sigma
_\varepsilon$). \citet{davis2012} provide numerical evidence that
the forecasting performance can be improved by using a log-likelihood
based loss function that incorporates information on the error
correlations. In this section, we consider both least squares and
log-likelihood estimates and study their theoretical properties.
A key contribution of our theoretical analysis is to verify suitable RE
and deviation conditions for the entire class of stable $\operatorname{VAR}(d)$ models.
Existing works either assume such conditions without verification, or
use a stringent condition on the model parameters, such as $\llVert  A\rrVert   < 1$,
as discussed in Section~\ref{secintro}.

We consider a single realization of $\{X^0, X^1, \ldots, X^T \}$
generated according to the VAR model (\ref{eqnVAR-defn}). We will
assume the error covariance matrix $\Sigma_{\varepsilon}$ is positive
definite so that $\Lambda_{\min} (\Sigma_{\varepsilon}) > 0$ and
$\Lambda_{\max}(\Sigma_\varepsilon) < \infty$.
We will also assume that the VAR process is \textit{stable}, that is,
$\det(\mathcal{A}(z)) \neq0$ on the unit circle $\{z \in\mathbb{C}\dvtx
\llvert  z\rrvert   = 1\}$. For stable $\operatorname{VAR}(d)$ processes, the spectral density (\ref
{eqnspectral-density-ARMA}) simplifies to
\[
\label{eqnspectral-density-AR} f_X (\theta) = \frac{1}{2\pi} \bigl(
\mathcal{A}^{-1} \bigl(e^{-i\theta
}\bigr) \bigr)
\Sigma_{\varepsilon} \bigl( \mathcal{A}^{-1} \bigl(e^{-i\theta}
\bigr) \bigr)^*.
\]
To deal with dependence in the VAR estimation problem, we will work
with $\mu_{\min}(\mathcal{A})$, $\mu_{\max}(\mathcal{A})$ and the
extreme eigenvalues of $\Sigma_\varepsilon$ instead of $\EuFrak
{m}(f_X)$ and $\mathcal{M}(f_X)$. For a $\operatorname{VAR}(d)$ process with serially
uncorrelated errors, equation (\ref{eqnmf-to-mu}) simplifies to
%
\begin{equation}
\label{eqnmf-to-mu-var} \mathcal{M}(f_X) \le\frac{1}{2\pi}
\frac{\Lambda_{\max}(\Sigma
_{\varepsilon})}{\mu_{\min} (\mathcal{A})},\qquad \EuFrak{m}(f_X) \ge\frac{1}{2\pi}
\frac{\Lambda_{\min}(\Sigma
_{\varepsilon})}{\mu_{\max} (\mathcal{A})}.
\end{equation}

\begin{figure}

\includegraphics{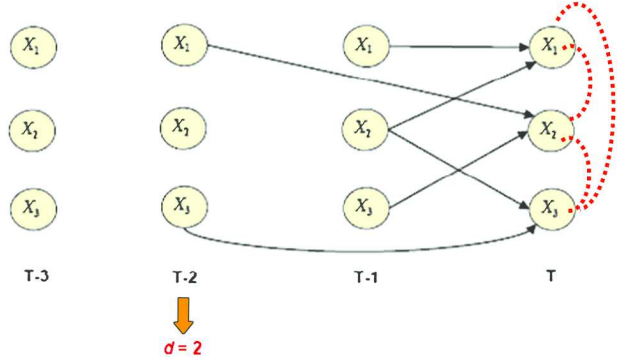}

\caption{Graphical representation of the VAR model (\protect\ref
{eqnVAR-defn}): directed edges (solid) correspond to the entries of the
transition matrices, undirected edges (dashed) correspond to the
entries of $\Sigma^{-1}_{\varepsilon}$.}\label{figl1-VAR-graph}
\end{figure}

This factorization helps provide better insight into the temporal and
contemporaneous dependence in VAR models. A graphical representation of
a stable $\operatorname{VAR}(d)$ model (\ref{eqnVAR-defn}) is provided in Figure~\ref
{figl1-VAR-graph}. The transition matrices $A_1, \ldots, A_d$ encode
the temporal dependence of the process. When the components of the
error process $\{ \varepsilon^t\}$ are correlated, $\Sigma_{\varepsilon
}^{-1}$ captures the additional contemporaneous dependence structure.
Expressing the estimation and prediction errors in terms of $\mu_{\min
}(\mathcal{A})$, $\mu_{\max}(\mathcal{A}), \Lambda_{\min}(\Sigma
_{\varepsilon})$ and $\Lambda_{\max}(\Sigma_{\varepsilon})$ instead of
$\EuFrak{m}(f_X)$ and $\mathcal{M}(f_X)$ help separate the effect of
the two sources of dependence.

%

We will often use the following alternative representation of a
$p$-dimensional $\operatorname{VAR}(d)$ process (\ref{eqnVAR-defn}) as a
$dp$-dimensional $\operatorname{VAR}(1)$ process $\tilde{X}^t = \tilde{A}_1   \tilde
{X}^{t-1} + \tilde{\varepsilon}^t$ with
%
\begin{eqnarray}\label{eqnvardto1}
\tilde{X}^t &=& \lleft[ \matrix{X^t
\cr
X^{t-1}
\cr
\vdots
\cr
X^{t-d+1}} \rright]_{dp \times1},\qquad
\tilde{A}_1 = \lleft[ \matrix{A_1 & A_2 &
\cdots& A_{d-1} & A_d
\cr
I_p & \mathbf{0} &
\cdots& \mathbf{0} & \mathbf{0}
\cr
\mathbf{0} & I_p & \cdots&
\mathbf{0} & \mathbf{0}
\cr
\vdots& \vdots& \ddots& \vdots& \vdots
\cr
\mathbf{0} &
\mathbf{0} & \cdots& I_p & \mathbf{0}} \rright]_{dp \times dp},\hspace*{-25pt}
\nonumber\\[-8pt]\\[-8pt]\nonumber
\tilde{\varepsilon}^t &=& \lleft[ \matrix{\varepsilon^t
\cr
\mathbf{0}
\cr
\vdots
\cr
\mathbf{0}} \rright]_{dp \times1}.\nonumber
\end{eqnarray}
The process $\tilde{X}^t$ with reverse characteristic polynomial
$\tilde{\mathcal{A}}(z):= I_{dp} - \tilde{A}_1 z$ is stable if and
only if the process $X^t$ is stable [\citet{lutkepohl2005new}].
However, the quantities $\mu_{\min}(\mathcal{A}),   \mu_{\max
}(\mathcal{A})$ are not necessarily the same as $\mu_{\min}(\tilde
{\mathcal{A}}),   \mu_{\max}(\tilde{\mathcal{A}})$.

\subsection{Estimation procedure}

Based on the data $\{X^0, \ldots, X^T \}$, we construct the following
regression problem:
\begin{eqnarray}
\label{eqnreg-setup} \underbrace{\lleft[ \matrix{\bigl(X^T
\bigr)'
\cr
\vdots
\cr
\bigl(X^{d}\bigr)'}
\rright]}_{\mathcal{Y}} & = & \underbrace{\lleft[ \matrix{
\bigl(X^{T-1}\bigr)' & \cdots& \bigl(X^{T-d}
\bigr)'
\cr
\vdots& \ddots& \vdots
\cr
\bigl(X^{d-1}
\bigr)' & \cdots& \bigl(X^0\bigr)'}
\rright]}_{\mathcal{X}} \underbrace{\lleft[ \matrix{A'_1
\cr
\vdots
\cr
A'_d} \rright]}_{B^*} +
\underbrace{\lleft[ \matrix{\bigl(\varepsilon^T\bigr)'
\cr
\vdots
\cr
\bigl(\varepsilon^d\bigr)'} \rright]}_{E},
\nonumber
\\
\operatorname{vec}(\mathcal{Y}) &=& \operatorname{vec}\bigl(\mathcal{X} B^*\bigr) + \operatorname{vec}(E),
\nonumber
\\
&=& (I \otimes\mathcal{X}) \operatorname{vec}\bigl(B^*\bigr) + \operatorname{vec}(E),
\nonumber
\\
\underbrace{Y}_{Np \times1} &=& \underbrace{Z}_{Np \times q} \underbrace{
\beta^*}_{q \times1} + \underbrace{\operatorname{vec}(E)}_{Np \times
1},\qquad N =
(T-d+1),\qquad q = dp^2,
\nonumber
\end{eqnarray}
with $N = T-d+1$ samples and $q = dp^2$ variables. We will assume that
$\beta^*$ is a $k$-sparse vector, that is, $\sum_{t=1}^d \llVert
\operatorname{vec}(A_t) \rrVert  _0 = k$.

We consider the following estimates for the transition matrices $A_1,
\ldots, A_d$, or equivalently, for $\beta^*$: (i) an
$\ell_1$-penalized least squares estimate of VAR coefficients ($\ell
_1$-LS), which does not exploit $\Sigma_{\varepsilon}$
%
\begin{equation}
\label{l1-LS} \mathop{\operatorname{argmin}}_{\beta\in\mathbb{R}^{q}} \frac
{1}{N}
\llVert Y - Z \beta\rrVert ^2 + \lambda_N \llVert \beta
\rrVert _1,
\end{equation}
%
and (ii) an $\ell_1$-penalized log-likelihood estimation ($\ell
_1$-LL) [\citet{davis2012}].
%
\begin{equation}
\label{l1-LL} \mathop{\operatorname{argmin}}_{\beta\in\mathbb{R}^{q}} \frac
{1}{N}
( Y - Z \beta )' \bigl( \Sigma_{\varepsilon}^{-1}
\otimes I \bigr) ( Y - Z \beta ) + \lambda_N \llVert \beta\rrVert
_1.
\end{equation}
This gives the maximum likelihood estimate of $\beta$, for
known $\Sigma_{\varepsilon}$. In practice, $\Sigma_\varepsilon$ is often
unknown and needs to be estimated from the data.
In the numerical experiments of Section~\ref{secsim}, we used the
residuals from a $\ell_1$-LS fit to estimate $\Sigma_\varepsilon$.
Further discussion on estimating $\Sigma_{\varepsilon}$ and a fast
algorithm based on block coordinate descent that minimizes (\ref
{l1-LL}) are presented in Appendix C (supplementary material [\citet{supp}]).

\subsection{Theoretical properties}\label{sectheory}
We analyze the estimates from optimization problems (\ref{l1-LS}) and
(\ref{l1-LL}) under a general penalized M-estimation framework
[\citet{powai2012}]. To motivate this general framework, note
that the VAR estimation problem with ordinary least squares is
equivalent to the following optimization:
%
\begin{equation}
\label{eqnLS} \mathop{\operatorname{argmin}}_{\beta\in\mathbb{R}^{q}} -2
\beta^{\prime} \hat{\gamma} + \beta^{\prime} \hat{\Gamma} \beta,
\end{equation}
where $\hat{\Gamma} =  ( I \otimes\mathcal{X}^{\prime} \mathcal{X}
/ N  )$, $\hat{\gamma} =  ( I \otimes\mathcal{X}^{\prime}
 ) Y/N $ are unbiased estimates for their population analogues.
A more general choice of $(\hat{\gamma}, \hat{\Gamma})$ in the
penalized version of the objective function leads to the following
optimization problem:
%
\begin{eqnarray}
\label{eqnpen-ML}
&\displaystyle\mathop{\operatorname{argmin}}_{\beta\in\mathbb{R}^{q}} -2
\beta^{\prime} \hat{\gamma} + \beta^{\prime} \hat{\Gamma} \beta+
\lambda_N \llVert \beta\rrVert _1,&
\nonumber\\[-8pt]\\[-8pt]\nonumber
&\displaystyle\hat{\Gamma} = \bigl(W \otimes\mathcal{X}^{\prime} \mathcal{X} / N
\bigr),\qquad \hat{\gamma} = \bigl( W \otimes\mathcal{X}^{\prime} \bigr) Y/N,&
\end{eqnarray}
where $W$ is a symmetric, positive definite matrix of weights.
Optimization problems (\ref{l1-LS}) and (\ref{l1-LL}) are special
cases of (\ref{eqnpen-ML}) with $W = I$ and $W = \Sigma
^{-1}_{\varepsilon}$, respectively.

First, we establish consistency of VAR estimates under the following
sufficient conditions: a modified restricted eigenvalue (RE)
[\citet{powai2012}] and a deviation condition. Then we show that
all stable VAR models satisfy these assumptions with high probability,
as long as the sample size is of the same order as required for consistency.



\begin{longlist}[(A1)]
\item[(A1)~Restricted eigenvalue (RE).] A symmetric matrix $\hat
{\Gamma}_{q \times q}$ satisfies restricted eigenvalue condition with
curvature $\alpha> 0$ and tolerance $\tau>0$ ($\hat{\Gamma} \sim
\operatorname{RE}(\alpha, \tau)$) if
%
\begin{equation}
\label{lower-RE} \theta{'} \hat{\Gamma} \theta\ge\alpha\llVert \theta
\rrVert ^2 - \tau\llVert \theta\rrVert ^2_1\qquad \forall \theta\in\mathbb{R}^q.
\end{equation}

The deviation condition ensures that $\hat{\gamma}$ and $\hat{\Gamma
}$ are well behaved in the sense that they concentrate nicely around
their population means. As $\hat{\gamma}$ and $\hat{\Gamma} \beta
^*$ have the same expectation, this assumption requires an upper bound
on their difference. Note that in the low-dimensional context of (\ref
{eqnLS}), $\hat{\gamma} - \hat{\Gamma} \beta^*$ is precisely
$\operatorname{vec}(\mathcal{X}'E)/N$.

\item[(A2)~Deviation condition.] There exists a deterministic
function $\mathbb{Q}  ( \beta^*, \Sigma_{\varepsilon}  )$
such that
%
\begin{equation}
\label{surr-dev} \bigl\llVert \hat{\gamma} - \hat{\Gamma} \beta^* \bigr\rrVert
_{\infty} \le\mathbb{Q} \bigl( \beta^*, \Sigma_{\varepsilon} \bigr) \sqrt
{\frac{\log d + 2 \log p}{N}}.
\end{equation}
\end{longlist}


%
\begin{prop}[(Estimation and prediction error)]\label{properror-bound}
Consider the penalized M-estimation problem (\ref{eqnpen-ML}) with $W
= I$ or $W = \Sigma_{\varepsilon}^{-1}$. Suppose $\hat{\Gamma}$
satisfies RE condition (\ref{lower-RE}) with $k \tau\le\alpha/ 32$,
and $(\hat{\Gamma}, \hat{\gamma})$ satisfies deviation bound (\ref
{surr-dev}).
Then, for any $\lambda_N \ge4 \mathbb{Q}(\beta^*, \Sigma_\varepsilon
) \sqrt{(\log d + 2 \log p)/ N}$, any solution $\hat{\beta}$ of
(\ref{eqnpen-ML}) satisfies
\begin{eqnarray*}
\bigl\llVert \hat{\beta} - \beta^* \bigr\rrVert _1 &\le&64 k
\lambda_N / \alpha,
\\
\bigl\llVert \hat{\beta} - \beta^* \bigr\rrVert &\le&16 \sqrt{k}
\lambda_N / \alpha,
\\
\bigl(\hat{\beta} - \beta^*\bigr)' \hat{\Gamma} \bigl(\hat{\beta}
- \beta^*\bigr) &\le&128 k \lambda_N^2 / \alpha.
\end{eqnarray*}
Further, a thresholded variant of lasso $\tilde{\beta} = \{\hat
{\beta}_j \mathbf{1}_{\llvert  \hat{\beta}_j\rrvert   > \lambda_N} \}
$ satisfies
\begin{eqnarray*}
&&\bigl\llvert \operatorname{supp}(\tilde{\beta}) \backslash \operatorname{supp}\bigl(\beta^*\bigr) \bigr\rrvert
\le \frac{192k}{\alpha_{\mathrm{RE}}}.
\end{eqnarray*}
\end{prop}

\begin{rems*}
(a) $\llVert   \hat{\beta} - \beta^* \rrVert  $ is precisely $\sum_{t = 1}^d \llVert
\hat{A}_t - A_t \rrVert  _F$, the $\ell_2$-error in estimating the
transition matrices. For $\ell_1$-LS, $(\hat{\beta} - \beta^*)'
\hat{\Gamma} (\hat{\beta} - \beta^*)$ is a measure of in-sample
prediction error under $\ell_2$-norm, defined by
$\sum_{t = d}^T \llVert  \sum_{h = 1}^d(\hat{A}_h - A_h) X^{t-h}\rrVert  ^2/N$.
For $\ell_1$-LL, $(\hat{\beta} - \beta^*)' \hat{\Gamma} (\hat
{\beta} - \beta^*) $ takes the form\break $\sum_{t = d}^T \llVert  \sum_{h =
1}^d (\hat{A}_h - A_h) X^{t-h}\rrVert  _{\Sigma_{\varepsilon}}^2/N$, where $\llVert
v\rrVert  _{\Sigma} := \sqrt{v' \Sigma^{-1} v}$. This can be viewed as a
measure of in-sample prediction error under a Mahalanobis-type distance
on $\mathbb{R}^p$ induced by~$\Sigma_{\varepsilon}$.

(b) The convergence rates are governed by two sets of parameters: (i)
dimensionality parameters, the dimension of the process $(p)$, order of
the process $(d)$, number of parameters $(k)$ in the transition
matrices $A_i$ and sample size $(N = T-d+1)$; (ii) internal parameters,
the curvature ($\alpha$), tolerance ($\tau$) and the deviation bound
$\mathbb{Q}(\beta^*, \Sigma_{\varepsilon})$. The squared $\ell
_2$-errors of estimation and prediction scale with the dimensionality
parameters as $k (2\log p + \log d)/N$, similar to the rates obtained
when the observations are independent [\citet
{bickel2009simultaneous}]. The temporal and cross-sectional dependence
affect the rates only through the internal parameters. Typically, the
rates are better when $\alpha$ is large and $\mathbb{Q}(\beta^*,
\Sigma_{\varepsilon}), \tau$ are small. In Propositions \ref
{verify-RE} and \ref{lemmartingaleconc}, we investigate in detail how
these quantities are related to the dependence structure of the process.

(c) Although the above proposition is derived under the assumption that
$d$ is the true order of the VAR process, the results hold even if $d$
is replaced by any upper bound $\bar{d}$ on the true order. This
follows from the fact that a $\operatorname{VAR}(d)$ model can also be viewed as
$\operatorname{VAR}(\bar{d})$, for any $\bar{d} > d$, with transition matrices $A_1,
\ldots, A_d, 0_{p \times p}, \ldots, {0}_{p \times p}$. Note that the
convergence rates change from $\sqrt{(\log p + 2\log d)/N}$ to $\sqrt
{(\log p + 2\log\bar{d})/N}$.

Proposition \ref{properror-bound} is deterministic; that is, it
assumes a fixed realization of $\{X^0, \ldots, X^T \}$. To show that
these error bounds hold with high probability, one needs to verify that
assumptions (A1--A2) are satisfied with high probability when $\{X^0,
\ldots, X^T \}$ is a random realization from the $\operatorname{VAR}(d)$ process. This
is accomplished in the next two propositions.
\end{rems*}

%
\begin{prop}[(Verifying RE for $\hat{\Gamma}$)]\label{verify-RE}
Consider a random realization $\{ X^0, \ldots, X^T\}$ generated
according to a stable $\operatorname{VAR}(d)$ process (\ref{eqnVAR-defn}). Then there
exist constants $c_i > 0$ such that for all $N \succsim\max\{\omega
^{2}, 1 \}k(\log d + \log p )$, 
with probability at least $1 - c_1 \exp(-c_2 N \min\{\omega^{-2}, 1
\})$, the matrix
\begin{eqnarray*}
\hat{\Gamma} = I_p \otimes\bigl(\mathcal{X}'
\mathcal{X}/N\bigr) \sim \operatorname{RE}(\alpha, \tau),
\end{eqnarray*}
where
\begin{eqnarray*}
\omega &=& c_3 \frac{\Lambda_{\max}(\Sigma_{\varepsilon})/\mu_{\min
}(\tilde{\mathcal{A}})}{\Lambda_{\min}(\Sigma_{\varepsilon})/\mu
_{\max}(\mathcal{A})},\qquad \alpha= \frac{\Lambda_{\min}(\Sigma_{\varepsilon})}{2 \mu_{\max
}(\mathcal{A})},
\\
\tau &=& \alpha\max\bigl\{\omega^{2}, 1 \bigr\}
\frac{\log d + \log p}{N}. 
\end{eqnarray*}
Further, if $\Sigma_{\varepsilon}^{-1}$ satisfies $\bar{\sigma
}^i_{\varepsilon} := \sigma^{ii}_{\varepsilon} - \sum_{j \neq i} \sigma
^{ij}_{\varepsilon} > 0$, for $i = 1, \ldots, p$, then, with the same
probability as above, the matrix
\[
\hat{\Gamma} = \Sigma_{\varepsilon}^{-1} \otimes \bigl(\mathcal
{X}'\mathcal{X}/N \bigr) \sim \operatorname{RE} \Bigl(\alpha \min
_{i} \bar {\sigma}^{i}_{\varepsilon}, \tau \max
_{i} \bar{\sigma }^i_{\varepsilon} \Bigr).
\]
\end{prop}

This proposition provides insight into the effect of temporal and
cross-sectional dependence on the convergence rates obtained in
Proposition \ref{properror-bound}. As mentioned earlier, the
convergence rates are faster for larger $\alpha$ and smaller $\tau$.
From the expressions of $\omega, \alpha$ and $\tau$, it is clear
that the VAR estimates have smaller error bounds when $\Lambda_{\max
}(\Sigma_{\varepsilon}), \mu_{\max}(\mathcal{A})$ are smaller and
$\Lambda_{\min}(\Sigma_{\varepsilon}), \mu_{\min}(\tilde{\mathcal
{A}})$ are larger, that is, when the spectrum is less spiky.

%
\begin{prop}[(Deviation bound)]\label{lemmartingaleconc} There exist
constants $c_i > 0$ such that for $N \succsim(\log d + 2 \log p)$,
with probability at least $1 - c_1 \exp [ -c_2 (\log d + 2\log
p) ]$, we have
\begin{eqnarray*}
&&\bigl\llVert \hat{\gamma} - \hat{\Gamma} \beta^* \bigr\rrVert _{\infty}
\le\mathbb{Q}\bigl(\beta^*, \Sigma_{\varepsilon}\bigr) \sqrt{\frac{\log d +
2\log p}{N}},
\end{eqnarray*}
where, for $\ell_1$-LS,
\[
\mathbb{Q}\bigl(\beta^*, \Sigma_{\varepsilon}\bigr) = c_0 \biggl[
\Lambda_{\max
}(\Sigma_{\varepsilon})+\frac{\Lambda_{\max}(\Sigma_{\varepsilon
})}{\mu_{\min} (\mathcal{A})}+
\frac{\Lambda_{\max}(\Sigma
_{\varepsilon})\mu_{\max}(\mathcal{A})}{\mu_{\min} (\mathcal{A})} \biggr]
\]
and for $\ell_1$-LL,
\[
\mathbb{Q}\bigl(\beta^*, \Sigma_{\varepsilon}\bigr) = c_0 \biggl[
\frac
{1}{\Lambda_{\min}(\Sigma_{\varepsilon})}+\frac{\Lambda_{\max
}(\Sigma_{\varepsilon})}{\mu_{\min} (\mathcal{A})}+\frac{\Lambda
_{\max}(\Sigma_{\varepsilon})\mu_{\max}(\mathcal{A})}{\Lambda_{\min
}(\Sigma_{\varepsilon}) \mu_{\min} (\mathcal{A})} \biggr].
\]
\end{prop}

As before, this proposition shows that the VAR estimates have lower
error bounds when $\Lambda_{\max}(\Sigma_{\varepsilon})$, $\mu_{\max
}(\mathcal{A})$ are smaller and $\Lambda_{\min}(\Sigma_{\varepsilon
})$, $\mu_{\min}(\mathcal{A})$ are larger, that is, when the
spectrum is less spiky.

\subsubsection*{Comparison with existing results}
The problem of sparse VAR
estimation has been theoretically studied in the literature in
[\citet{songbickel2011,Chudik2011IVAR,WuWu2014}]. Next, we
briefly highlight differences between our results and these works.
First, the results of \citet{Chudik2011IVAR} rely on a~priori
available neighborhood information for every time series, {which
implies that the structure of transition matrices $\{A_t\}_{t=1}^d$ is
known, and only their magnitudes need to be estimated}. This is a
significant limitation compared to regularized methods like lasso,
which do not require any prior knowledge on the sparsity pattern in the
transition matrices. The theoretical upper bounds on VAR estimation
error established in \citet{songbickel2011} do not decrease as
the sample size $T$ increases, and hence do not ensure consistency
beyond very strict conditions. Also, the results in their paper and in
\citet{WuWu2014} are established assuming RE holds, while a
significant portion of our analysis is devoted to establish that RE and\vadjust{\goodbreak}
deviation bounds hold with high probability. We also provide in-depth
analysis on how the relevant constants are affected by the dependence
present in the data. Finally, our work is the first one to provide
theoretical analysis of the log-likelihood based VAR estimation
procedure, which does not fit directly into the regression setting
considered in the aforementioned papers.

\section{Extension to other regularized estimation problems}\label{secextension}
The deviation inequalities established in Section~\ref
{secmain-results} can be easily integrated with the vast body of
existing literature of high-dimensional statistics for i.i.d. data and
study other regularized estimation problems in the context of
high-\break dimensional time series. To demonstrate this, in this section we
establish consistency of sparse covariance estimation by
hard-thresholding [\citet{BL08covthresh}] for high-dimensional
time series and discusss the main steps in extending the results to
some nonconvex penalties for sparse regression and group lasso and
nuclear norm penalties for inducing structured sparsity.

\subsection{Sparse covariance estimation}\label{seccov-est}
Consider a $p$-dimensional centered Gaussian stationary time series $\{
X^t \}_{t \in\mathbb{Z}}$ satisfying Assumption \ref
{assumpspectral-density}. Based on realizations $\{X^1, \ldots, X^n \}
$ generated according to the above stationary process, we aim to
estimate the contemporaneous covariance matrix $\Sigma= \Gamma(0)$.
The sample covariance matrix $\hat{\Gamma}(0) = \frac{1}{n} \sum_{t=1}^n (X^t - \bar{X}) (X^t - \bar{X})' $ is known to be
inconsistent when $p$ grows faster than $n$.
\citet{BL08covthresh} showed that when the samples are generated
independently from a centered Gaussian or subGaussian distribution, a
thresholded version of the sample covariance matrix $T_u(\hat{\Gamma
}(0)) = \{ \hat{\Gamma}_{ij}(0) \mathbf{1}_{\llvert  \hat{\Gamma}_{ij}(0)\rrvert
> u} \}$ can perform consistent estimation if ${\Gamma(0)}$ belongs to
the following uniformity class of approximately sparse matrices:
\[
\label{defncov-class} \mathcal{U}_{\tau}\bigl(q, c_0(p), M\bigr):=
\Biggl\{ \Sigma\dvtx  \sigma_{ii} \le M, \sum_{j=1}^p
\llvert \sigma_{ij}\rrvert ^q \le c_0(p),
\mbox{ for all } i \Biggr\}.
\]
Next, we establish consistent estimation for time series data, provided
that the underlying process is stable. The effect of dependence on the
estimation accuracy is captured by the stability measures introduced in
Section~\ref{secmain-results}. Asymptotic theory for sparse covariance
estimation was also considered in \citet{chen-wu-13-AOS},
assuming a decay on the functional dependence measure.

\begin{prop}\label{propcov-est}
Let $\{X^t\}_{t=1}^n$ be generated according to a $p$-\break dimensional
stationary centered Gaussian process with spectral density $f_X$,
satisfying Assumption \ref{assumpspectral-density}. Then, uniformly on
$\mathcal{U}_{\tau} (q, c_0(p), M)$, for sufficiently large $M'$, if
$u_n =\mathcal{M}(f_X, 2) M' \sqrt{\log p/n}$ and
$n \succsim\mathcal{M}^2(f_X, 2) \log p $, then
\begin{eqnarray*}
\bigl\llVert T_{u_n} \bigl(\hat{\Gamma}(0)\bigr) - \Gamma(0) \bigr
\rrVert &=& O_p \biggl(c_0(p) \biggl(\mathcal{M}^2(f_X,
2) \frac{ \log p}{n} \biggr)^{\vfrac
{1-q}{2}} \biggr),
\\
\frac{1}{p}\bigl\llVert T_{u_n} \bigl(\hat{\Gamma}(0)\bigr) -
\Gamma(0) \bigr\rrVert _F &=& O_p \biggl(c_0(p)
\biggl( \mathcal{M}^2(f_X, 2) \frac{ \log
p}{n}
\biggr)^{1-\sklfrac{q}{2}} \biggr).
\end{eqnarray*}
\end{prop}

\subsection{Sparse regression with nonconvex penalties}\label{subsecnonconvex-reg}
There is a vast body of literature on regularized regression using
nonconvex penalties for i.i.d. data [\citet{fan2001scad,zhang2010mcp}]. A recent line work has derived unified theoretical
treatments of these procedures and compared their estimation accuracy
to convex procedures such as lasso [\citet{fanlv2013eqv,lowai2013}]. These results indicate that in certain high-dimensional
regimes, the estimation error of nonconvex penalties like SCAD, MCP
scales roughly in the same order as lasso. Next, we argue that similar
conclusions hold for time series models, as well.

Consider a stochastic regression problem of Section~\ref{secstoch-reg}
subject to a SCAD or MCP penalty. \citet{lowai2013} establish
that under suitable restricted strong convexity (RSC) condition on the
loss function $\mathcal{L}_n(\cdot)$, if the sup norm of the gradient $\llVert
\nabla\mathcal(L)_n(\beta^*) \rrVert  _{\infty}$ scales with $\sqrt{\log
p/n}$, then any local solution of the penalized objective function has
an estimation error at most $O(\sqrt{k \log p/n})$. For the choice of
a least squares loss function, $\mathcal{L}_n(\beta)=\llVert  \mathcal
{Y}-\mathcal{X}\beta\rrVert  ^2/2n$ and $\nabla\mathcal{L}_n(\beta^*) =
-\mathcal{X}'E/n$.

Since the loss function is convex, their RSC takes the form
\[
\frac{1}{n} \frac{\llVert  \mathcal{X} \Delta\rrVert  ^2}{\llVert  \Delta\rrVert  ^2} \ge \alpha_1 \llVert \Delta
\rrVert ^2 - \tau_1 \frac{\log p}{n} \llVert \Delta
\rrVert _1^2\qquad\mbox{for all } \llVert \Delta\rrVert
\le1.
\]
This is in the spirit of the RE conditions verified in Section~\ref
{secVAR} and can be proven using similar discretization arguments
presented in this paper, if we assume $\Gamma(0)$ satisfies an RE
condition with the restricted eigenvalue $\alpha_1$ is at least as
large as $1/(a-1)$ for SCAD and $1/b$ for MCP.

The deviation condition on $\llVert  \nabla\mathcal{L}_n(\beta^*) \rrVert
_{\infty}$ is identical to the one considered in this paper, and the
results presented here are directly applicable.
\subsection{Regularized regression with structured sparsity}\label
{subsecgrp-lowrank}
In a recent review paper, \citet{negahban2012unified} established
a unified framework to analyze a class of decomposable penalties. This
includes the popular group lasso penalty for high-dimensional
regression under structured sparsity and nuclear norm penalty for
matrix estimation under low-rank assumption. In a time series context,
these methods have been proposed in the literature to incorporate
information on different economic sectors and the assumption of latent
factors driving the market [\citet{songbickel2011,negwai2011}].
As before, the theoretical results rely crucially on two key
conditions: a restricted strong convexity on the loss function and a
suitable deviation bound on the gradient. The restricted eigenvalue
assumption for group lasso can be verified using the deviation
inequalities of Proposition \ref{propconc-S} and a discretization
argument modified for group structures. The deviation inequalities can
be derived along the same line. For low-rank modeling of $\operatorname{VAR}(1)$
process, we can prove that the minimum eigenvalue of $\mathcal
{X}'\mathcal{X}/N$ is bounded away from zero with high probability,
and the deviation bounds on the\vadjust{\goodbreak} operator norm of $\mathcal{X}'E/N$ can
be established using the deviation inequality of~(\ref
{eqnconc-XY-var}) and a discretization argument presented in
[\citet{mythesis}]. This leads to new results on group lasso for
stochastic regression and extends the results of \citet
{negwai2011} to the entire class of stable $\operatorname{VAR}(1)$ models. We leave the
details to the reader, as the proofs follow the same road map used in
this paper.

\section{Numerical experiments}\label{secsim}
\subsection{Stochastic regression}

\begin{figure}

\includegraphics{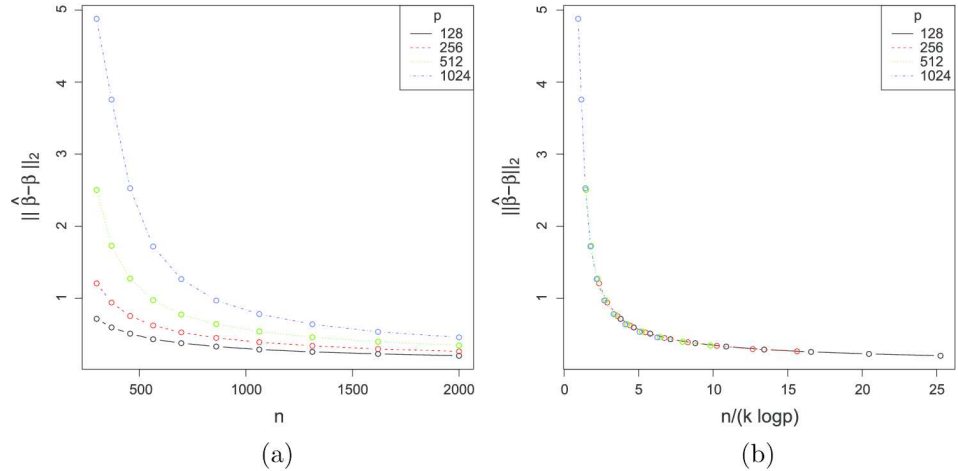}

\caption{Estimation error of lasso $\llVert  \hat{\beta} -\beta^*\rrVert  $ in
stochastic regression with serially correlated error.
Note that the error curves align perfectly, showing the errors scale as
$\sqrt{k \log p / n}$.
\textup{(a)}~$\llVert   \hat{\beta} - \beta^*\rrVert  $ vs. $n$,
\textup{(b)}~$\llVert   \hat{\beta} - \beta^*\rrVert  $ vs. $n / k \log p$.}\label{figstoch-reg-scaling}
\end{figure}

In this experiment, we demonstrate how the estimation error of lasso
scales with $n$ and $p$, when the dependence parameters do not change.
We simulate predictors from a $p$-dimensional
($p = 128, 264, 512, 1024$) stationary process $\{X^t\}$ with\vspace*{1pt}
independent components following a Gaussian $\operatorname{AR}(2)$ process $X^t_i =
1.2X^{t-1}_i - 0.36 X^{t-2}_i + \xi^t$, $\Gamma_{X_j}(0)=1$. We
simulate the errors $\{ \varepsilon^t\}$ according to a univariate
MA($2$) process $\varepsilon^t = \eta^t - 0.8 \eta^{t-1} + 0.16 \eta
^{t-2}$, $\{\eta^t\}$ Gaussian white noise. For different values of
$p$, we generate sparse coefficient vectors $\beta^*$ with $k \approx
\sqrt{p}$ nonzero entries, with a signal-to-noise ratio of $1.2$.
Using a tuning parameter $\lambda_n = \sqrt{\log p / n}$, we apply
lasso on simulated samples of size $n \in(100, 3000)$. The $\ell
_2$-error of estimation $\llVert  \hat{\beta} - \beta^* \rrVert  $ is depicted in
Figure~\ref{figstoch-reg-scaling}. The left panel displays the errors for
different values of $p$, plotted against the sample size~$n$. As
expected, the errors are larger for larger $p$. The right panel
displays the estimation errors against the rescaled sample size $n/k
\log p$. The error curves for different values of $p$ now align very
well. This demonstrates that lasso can achieve an estimation error rate
of $\sqrt{k \log p / n}$, even with stochastic predictors and serially
correlated errors.

%
\subsection{VAR estimation}
We evaluate the performance of $\ell_1$-LS and $\ell_1$-LL on
simulated data and compare it with the performance of ordinary least
squares (OLS) and Ridge estimates. Implementing $\ell_1$-LL requires
an estimate of $\Sigma_{\varepsilon}$ in the first step. We use the
residuals from $\ell_1$-LS to construct a plug-in estimate $\hat
{\Sigma_{\varepsilon}}$. To evaluate the effect of error correlation on
the transition matrix estimates more precisely, we also implement an
oracle version, $\ell_1$-LL-O, which uses the true $\Sigma_{\varepsilon
}$ in the estimation. Next, we describe the simulation settings, choice
of performance metrics and discuss the results.

We design two sets of numerical experiments: (a) SMALL VAR ($p = 10, d
= 1, T = 30, 50$) and (b) MEDIUM VAR ($p = 30, d = 1, T = 80, 120, 160$).
In each setting, we generate an adjacency matrix $A_1$ with $5\sim10\%
$ nonzero edges selected at random and rescale to ensure that the
process is stable with $SNR = 2$. We generate three different error
processes with covariance matrix $\Sigma_{\varepsilon}$ from one of the
following families:
\begin{longlist}[(2)]
\item[(1)] Block-I: $\Sigma_\varepsilon= ((\sigma_{\varepsilon, ij}))_{1
\le i, j \le p}$ with $\sigma_{\varepsilon, ii} = 1$, $\sigma_{\varepsilon
, ij} = \rho$ if $1 \le i \neq j \le p/2$, $\sigma_{\varepsilon, ij} =
0$ otherwise;
\item[(2)] Block-II: $\Sigma_\varepsilon= ((\sigma_{\varepsilon,
ij}))_{1 \le i, j \le p}$ with $\sigma_{\varepsilon, ii} = 1$, $\sigma
_{\varepsilon, ij} = \rho$ if $1 \le i \neq j \le p/2$ or $p/2 < i \neq
j \le p$, $\sigma_{\varepsilon, ij} = 0$ otherwise;
\item[(3)] Toeplitz: $\Sigma_\varepsilon= ((\sigma_{\varepsilon,
ij}))_{1 \le i, j \le p}$ with $\sigma_{\varepsilon, ij} = \rho^{\llvert  i-j\rrvert  }$.
\end{longlist}
%
%
\begin{figure}

\includegraphics{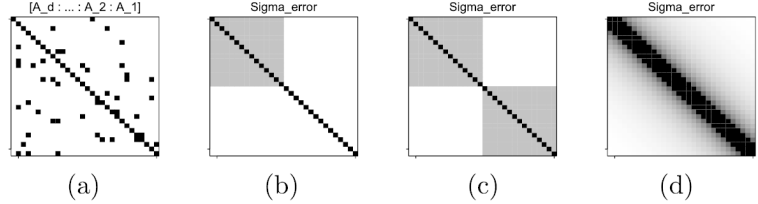}

\caption{Adjacency matrix $A_1$ and error covariance matrix $\Sigma
_\varepsilon$ of different types used in the simulation studies.
\textup{(a)}~$A_1$,
\textup{(b)} $\Sigma_\varepsilon$: Block-\textup{I},
\textup{(c)} $\Sigma_\varepsilon$: Block-\textup{II},
\textup{(d)} $\Sigma_\varepsilon$: Toeplitz.}\label{figsimsetting}
\end{figure}
We let $\rho$ vary in $\{0.5, 0.7, 0.9 \}$. Larger values of $\rho$
indicate that the error processes are more strongly correlated.
Figure~\ref{figsimsetting} illustrates the structure of a random
transition matrix used in our simulation and the three different types
of error covariance structures.

We compare the different methods for VAR estimation (OLS, $\ell_1$-LS,
$\ell_1$-LL, $\ell_1$-LL-O, Ridge) based on the following performance metrics:
\begin{longlist}[(2)]
\item[(1)] \textit{Model Selection}. Area under ROC curve (AUROC);
\item[(2)] \textit{Estimation error}. Relative estimation accuracy
$\llVert   \hat{A}_1 - A_1 \rrVert  _F / \llVert   A_1 \rrVert  _F$.
\end{longlist}

We report the results for small VAR with $T =30$ and medium VAR with $T
= 120$ averaged over $1000$ replicates in Tables~\ref{tblvarsmall} and
\ref{tblvarmed}. The results in the other settings are qualitatively
similar, although the overall accuracy changes with the sample size. We
find that the regularized VAR estimates outperform ordinary least
squares uniformly in all the cases.

In terms of model selection, the $\ell_1$-penalized estimates perform
fairly well, as reflected in their AUROC. OLS and ridge regression do
not perform any model selection. Further, for all three choices of
$\Sigma_\varepsilon$, the two variants of $\ell_1$-LL outperform $\ell
_1$-LS. The difference in their performance is more prominent for
larger values of~$\rho$. Among the three covariance structures, the
difference between LS- and \mbox{LL-}based methods is more prominent in the
Block-II and Toeplitz families, since the error processes are more
strongly correlated. Finally, in all cases, the accuracy of $\ell
_1$-LL lies between $\ell_1$-LS and $\ell_1$-LL-O, which suggests
that a more accurate estimation of $\Sigma_{\varepsilon}$ might improve
the model selection performance of regularized VAR estimates.

In terms of estimation error, the conclusions are broadly the same. The
effect of over-fitting is reflected in the performance of ordinary
least squares. In many\vadjust{\goodbreak} settings, the estimation error of ordinary least
squares is even twice as large as the signal strength. The performance
of ordinary least squares deteriorates when the error processes are
more strongly correlated; see, for example, $\rho= 0.9$ for block-II.
Ridge regression performs better than ordinary least squares, as it
applies shrinkage on the coefficients. However, the $\ell_1$-penalized
estimates show higher accuracy than Ridge in almost all cases. This is
somewhat expected as the data were simulated from a sparse model with
strong signals, whereas Ridge regression tends to favor a nonsparse
model with many small coefficients.



\begin{table}
\tabcolsep=0pt
\caption{$\operatorname{VAR}(1)$ model with $p = 10$, $T = 30$}\label{tblvarsmall}
\begin{tabular*}{\tablewidth}{@{\extracolsep{\fill}}@{}lcccccccccc@{}}
\hline
& & \multicolumn{3}{c}{\textbf{Block-I}} & \multicolumn{3}{c}{\textbf{Block-II}} & \multicolumn{3}{c@{}}{\textbf{Toeplitz}}\\[-6pt]
& & \multicolumn{3}{c}{\hrulefill} & \multicolumn{3}{c}{\hrulefill} & \multicolumn{3}{c@{}}{\hrulefill}\\
&$\bolds{\rho}$ & $\bolds{0.5}$ & $\bolds{0.7}$ & $\bolds{0.9}$ & $\bolds{0.5}$ & $\bolds{0.7}$ & $\bolds{0.9}$ & $\bolds{0.5}$ & $\bolds{0.7}$ & $\bolds{0.9}$\\
\hline
AUROC & $\ell_1$-LS & 0.78 & 0.77 & 0.74 & 0.74 & 0.7\phantom{0} & 0.64 & 0.76 & 0.72 & 0.63\\
& $\ell_1$-LL & 0.79 & 0.79 & 0.76 & 0.77 & 0.77 & 0.76 & 0.78 & 0.76 & 0.74 \\
& $\ell_1$-LL-O & 0.84 & 0.83 & 0.8\phantom{0} & 0.82 & 0.82 & 0.82 & 0.83 & 0.82 & 0.8\phantom{0}
\\[3pt]
Estimation & OLS & 1.51 & 1.67 & 2.31 & 1.73 & 2.16 & 3.57 & 1.7\phantom{0} & 2.14 & 3.57 \\
\quad error & $\ell_1$-LS & 0.74 & 0.75 & 0.76 & 0.77 & 0.8\phantom{0} & 0.87 & 0.77 & 0.8\phantom{0} & 0.88 \\
& $\ell_1$-LL & 0.7\phantom{0} & 0.7\phantom{0} & 0.69 & 0.73 & 0.72 & 0.72 & 0.73 & 0.73 & 0.74 \\
& $\ell_1$-LL-O & 0.65 & 0.64 & 0.63 & 0.66 & 0.65 & 0.63 & 0.66 & 0.66 & 0.65 \\
& Ridge & 0.78 & 0.78 & 0.79 & 0.77 & 0.78 & 0.8\phantom{0} & 0.8\phantom{0} & 0.82 & 0.85 \\
\hline
\end{tabular*}
\end{table}

%
\begin{table}[b]
\tabcolsep=0pt
\caption{$\operatorname{VAR}(1)$ model with $p = 30$, $T = 120$}\label{tblvarmed}
\begin{tabular*}{\tablewidth}{@{\extracolsep{\fill}}@{}lcccccccccc@{}}
\hline
& & \multicolumn{3}{c}{\textbf{Block-I}} & \multicolumn{3}{c}{\textbf{Block-II}} & \multicolumn{3}{c@{}}{\textbf{Toeplitz}}\\[-6pt]
& & \multicolumn{3}{c}{\hrulefill} & \multicolumn{3}{c}{\hrulefill} & \multicolumn{3}{c@{}}{\hrulefill}\\
&$\bolds{\rho}$ & $\bolds{0.5}$ & $\bolds{0.7}$ & $\bolds{0.9}$ & $\bolds{0.5}$ & $\bolds{0.7}$ & $\bolds{0.9}$ & $\bolds{0.5}$ & $\bolds{0.7}$ & $\bolds{0.9}$\\
\hline
AUROC & $\ell_1$-LS & 0.91 & 0.87 & 0.8\phantom{0} & 0.82 & 0.75 & 0.63 & 0.92 & 0.88 & 0.77 \\
& $\ell_1$-LL & 0.91 & 0.89 & 0.85 & 0.85 & 0.85 & 0.85 & 0.93 & 0.92 & 0.91 \\
& $\ell_1$-LL-O & 0.93 & 0.91 & 0.87 & 0.88 & 0.88 & 0.88 & 0.95 & 0.94 & 0.92
\\[3pt]
Estimation & OLS & 1.65 & 1.91 & 2.74 & 2.33 & 2.98 & 4.94 & 1.77 & 2.24 & 3.74 \\
\quad error & $\ell_1$-LS & 0.68 & 0.73 & 0.8\phantom{0} & 0.83 & 0.9\phantom{0} & 0.98 & 0.68 & 0.72 & 0.85 \\
& $\ell_1$-LL & 0.67 & 0.67 & 0.67 & 0.78 & 0.77 & 0.74 & 0.65 & 0.62 & 0.57 \\
& $\ell_1$-LL-O & 0.63 & 0.63 & 0.63 & 0.74 & 0.73 & 0.7\phantom{0} & 0.61 & 0.57 & 0.52 \\
& Ridge & 0.8\phantom{0} & 0.81 & 0.83 & 0.86 & 0.89 & 0.92 & 0.8\phantom{0} & 0.82 & 0.86 \\
\hline
\end{tabular*}\vspace*{-6pt}
\end{table}

\section{Discussion}\label{secdiscussion}
In this paper, we consider the theoretical properties of regularized
estimates in sparse high-dimensional time series models when the data
are generated
from a multivariate stationary Gaussian process. The Gaussian
assumption could be conceived as a limiting factor, since interesting
models including
regression with categorical predictors, VAR estimation with
heavy-tailed and/or heteroscedastic errors, and popular models
exhibiting nonlinear
dependences such as ARCH and GARCH are not covered. Note, however, that
the only place in the analysis where the Gaussian assumption is used
is in developing the concentration bound of $S$ around its expectation
$\Gamma(0)$. Since the spectral density characterizes the entire
distribution for this class, it has direct implications on the
concentration behavior. For nonlinear and/or non-Gaussian processes,
one needs to control higher order dependence, and changing to higher
order spectra could potentially be useful. Although the use of
covariance and higher order spectra is common in developing limit
theorems of low-dimensional stationary process [\citet
{rosenblatt1985stationary,giraitis2012large}], developing a suitable
concentration bound for nonlinear/non-Gaussian dependence designs is
not a trivial problem and is left as a key topic for future developments.

\section*{Acknowledgements}
We thank the Editor Runze Li, the Associate Editor and three anonymous
reviewers, whose comments led to several improvements in the paper.

%
\begin{supplement}[id=suppA]
\stitle{Supplement to ``Regularized estimation in sparse high-dimensional time series models''}
\slink[doi]{10.1214/15-AOS1315SUPP} 
\sdatatype{.pdf}
\sfilename{aos1315\_supp.pdf}
\sdescription{For the sake of brevity, we moved the appendices
containing many of the technical proofs and detailed discussions to the
supplementary document [\citet{supp}].}
\end{supplement}





%

\printaddresses
\end{document}